\documentclass[a4paper, 11pt, twoside, final]{scrartcl}

\usepackage{calc}
\usepackage{ifthen}
\usepackage{xargs}

\usepackage[T1]{fontenc}
\usepackage[utf8]{inputenc}
\usepackage[ngerman,british]{babel}

\usepackage[intlimits]{amsmath}
\usepackage{amssymb}
\usepackage{mathtools}
\usepackage[framed, hyperref, thref, amsmath]{ntheorem}
\usepackage[svgnames]{xcolor}
\usepackage[ntheorem, framemethod=tikz]{mdframed}
\usepackage{units}
\usepackage{mleftright}
\allowdisplaybreaks

\usepackage[libertine,slantedGreek]{newtxmath} 
\usepackage[bb=px,bbscaled=0.94,cal=dutchcal,calscaled=0.95]{mathalfa}
\usepackage{bm}
\usepackage[osf,p]{libertine}

\usepackage[headsepline,footsepline,plainfootsepline,
headsepline=1pt,footsepline=1pt]{scrlayer-scrpage}
\clearscrheadfoot
\ohead{\pagemark}
\cehead{\footnotesize\itshape\scshape \SHORTAUTHORsave}
\cohead{\footnotesize\itshape\scshape \SHORTTITLEsave}

\usepackage[shortlabels]{enumitem}
\setlist{listparindent=\parindent,parsep=\parskip}

\usepackage{graphicx}
\usepackage{placeins}
\usepackage{caption}
\usepackage{subfig}
\DeclareCaptionFont{gray}{\color{black}}
\usepackage{array}
\usepackage{colortbl}
\usepackage{longtable}
\usepackage{ltablex}
\usepackage{varioref}
\usepackage{hyperref}
\usepackage{breakurl}
\usepackage{bookmark}
\hypersetup{hidelinks}
\hypersetup{breaklinks=true}
\hypersetup{bookmarksopen,bookmarksopenlevel=1}
\addto\extrasbritish{%
}

\KOMAoption{BCOR}{5mm}
\KOMAoption{DIV}{classic}
\KOMAoption{headinclude}{false}
\KOMAoption{footinclude}{false}
\KOMAoption{pagesize}{auto}
\recalctypearea

\KOMAoption{headings}{small}
\KOMAoption{numbers}{autoendperiod}
\addtokomafont{section}{\scshape\liningnums}
\addtokomafont{subsection}{\scshape\liningnums}
\addtokomafont{subsubsection}{\scshape\liningnums}

\RedeclareSectionCommand[afterskip=-1em]{subsection}
\RedeclareSectionCommand[afterskip=-1em]{subsubsection}

\newcommand{\AUTHORsave}{}
\newcommand{\AUTHOR}[1]{\renewcommand{\AUTHORsave}{#1}}
\newcommand{\SHORTAUTHORsave}{}
\newcommand{\SHORTAUTHOR}[1]{%
  \ifthenelse{\equal{#1}{}}%
  {\renewcommand{\SHORTAUTHORsave}{\AUTHORsave}}%
  {\renewcommand{\SHORTAUTHORsave}{#1}}%
}
\newcommand{\TITLEsave}{}
\newcommand{\TITLE}[1]{\renewcommand{\TITLEsave}{#1}}
\newcommand{\SHORTTITLEsave}{}
\newcommand{\SHORTTITLE}[1]{%
  \ifthenelse{\equal{#1}{}}%
  {\renewcommand{\SHORTTITLEsave}{\TITLEsave}}%
  {\renewcommand{\SHORTTITLEsave}{#1}}%
}
\newcommand{\INSTITUTEsave}{}
\newcommand{\INSTITUTE}[1]{\renewcommand{\INSTITUTEsave}{#1}}
\newcommand{\CORRESPONDENCEsave}{}
\newcommand{\CORRESPONDENCE}[1]{\renewcommand{\CORRESPONDENCEsave}{#1}}
\newcommand{\ABSTRACTsave}{}
\newcommand{\ABSTRACT}[1]{\renewcommand{\ABSTRACTsave}{#1}}
\newcommand{\KEYWORDSsave}{}
\newcommand{\KEYWORDS}[1]{\renewcommand{\KEYWORDSsave}{#1}}
\newcommand{\AMSCLASSsave}{}
\newcommand{\AMSCLASS}[1]{\renewcommand{\AMSCLASSsave}{#1}}
\newcommand{\Inst}[1]{\textsuperscript{#1}}
\newcommand{\email}[1]{\href{mailto:#1}{\nolinkurl{#1}}}

\newcommand{\TitleHeader}{%
  \noindent%
  \rule{\linewidth}{1pt}\\[10pt]%
  {\sffamily\bfseries\Large\TITLEsave}\\[10pt]%
  {\sffamily\bfseries\small\AUTHORsave}\\%
  {\sffamily\scriptsize\INSTITUTEsave}\\[15pt]%
  \parbox{\linewidth}{\sffamily\scriptsize {\bfseries
      Correspondence}\\\CORRESPONDENCEsave}\\[10pt]%
  \rule{\linewidth}{1pt}\\[15pt]%
  \hspace*{40pt}%
  \parbox{\linewidth-40pt}{%
    \small%
    {\sffamily\bfseries Abstract.} \quad \ABSTRACTsave%
  }\\[5pt]%
  \hspace*{40pt}%
  \parbox{\linewidth-40pt}{%
    \small%
    {\sffamily\bfseries Keywords.} \quad \KEYWORDSsave%
  }\\[5pt]%
  \hspace*{40pt}%
  \parbox{\linewidth-40pt}{%
    \small%
    {\sffamily\bfseries AMS subject classification.} \quad \AMSCLASSsave%
  }\\[5pt]%
  \rule{\linewidth}{1pt}%
}

\newcommand{\qed}{}

\newlength{\frametopsep}
\definecolor{ThmFrame}{gray}{0.95}
\definecolor{ThmBody}{gray}{0.95}
\mdfdefinestyle{thmframe}{innerleftmargin=5pt, innerrightmargin=5pt,
  innermargin=-5.5pt, outermargin=-5.5pt,
  roundcorner=5pt, linewidth=0.5pt, linecolor=ThmFrame,
  backgroundcolor=ThmBody, splittopskip=12.5pt, splitbottomskip=5pt}
\makeatletter
\newtheoremstyle{plainheaderbreak}
  {\item[]{\theorem@headerfont ##1\ ##2\theorem@separator}\hskip\labelsep}%
  {\item[]{\theorem@headerfont ##1\ ##2\
    (##3)}\theorem@separator\hskip\labelsep}
\makeatother

\theoremheaderfont{\sffamily\bfseries\scshape}
\theorembodyfont{}
\theoremstyle{plainheaderbreak}
\theoremseparator{.}
\theoremframepreskip{\frametopsep}
\theoremframepostskip{\frametopsep}
\newmdtheoremenv[style=thmframe]{Definition}{Definition}[section]
\newmdtheoremenv[style=thmframe]{Problem}[Definition]{Problem}
\newmdtheoremenv[style=thmframe]{Assumption}[Definition]{Assumption}
\newmdtheoremenv[style=thmframe]{Annahme}[Definition]{Annahme}

\theoremheaderfont{\sffamily\bfseries\scshape\upshape}
\theorembodyfont{\itshape}

\newmdtheoremenv[style=thmframe]{Theorem}[Definition]{Theorem}
\newmdtheoremenv[style=thmframe]{Satz}[Definition]{Satz}
\newmdtheoremenv[style=thmframe]{Proposition}[Definition]{Proposition}
\newmdtheoremenv[style=thmframe]{Lemma}[Definition]{Lemma}
\newmdtheoremenv[style=thmframe]{Corollary}[Definition]{Corollary}
\newmdtheoremenv[style=thmframe]{Korollar}[Definition]{Korollar}

\theoremheaderfont{\sffamily\bfseries\scshape\upshape}
\theorembodyfont{\sffamily}
\theoremstyle{break}

\newmdtheoremenv[style=thmframe]{Algorithm}[Definition]{Algorithm}
\newmdtheoremenv[style=thmframe]{Algorithmus}[Definition]{Algorithmus}

\newlength{\unframetopsep}
\theoremheaderfont{\itshape}
\theorembodyfont{}
\theoremstyle{plainheaderbreak}
\theorempreskip{\unframetopsep}
\theorempostskip{\unframetopsep}

\theoremprework{%
  \renewcommand{\qed}{%
    \hspace*{\fill}\nolinebreak[3]%
    \nopagebreak[3]\hspace*{\fill}{\raisebox{0.25pt}{\scalebox{0.9}{$\medcirc$}}}}}
\theorempostwork{}

\theoremprework{%
  \renewcommand{\qed}{%
    \hspace*{\fill}\nolinebreak[3]%
    \nopagebreak[3]\hspace*{\fill}{\raisebox{0.25pt}{\scalebox{0.9}{$\medcirc$}}}}}
\theorempostwork{}

\theoremprework{%
  \renewcommand{\qed}{%
    \hspace*{\fill}\nolinebreak[3]%
    \nopagebreak[3]\hspace*{\fill}{\raisebox{0.25pt}{\scalebox{0.9}{$\medcirc$}}}}}
\theorempostwork{}

\newtheorem{Counterexample}[Definition]{Counterexample}

\theoremprework{%
  \renewcommand{\qed}{%
    \hspace*{\fill}\nolinebreak[3]%
    \nopagebreak[3]\hspace*{\fill}{\raisebox{0.25pt}{\scalebox{0.9}{$\medcirc$}}}}}
\theorempostwork{}

\theoremprework{%
  \renewcommand{\qed}{%
    \hspace*{\fill}\nolinebreak[3]%
    \nopagebreak[3]\hspace*{\fill}{\raisebox{0.25pt}{\scalebox{0.9}{$\medcirc$}}}}}
\theorempostwork{}

\newtheorem{Remark}[Definition]{Remark}

\theoremprework{%
  \renewcommand{\qed}{%
    \hspace*{\fill}\nolinebreak[3]%
    \nopagebreak[3]\hspace*{\fill}{\raisebox{0.25pt}{\scalebox{0.9}{$\medcirc$}}}}}
\theorempostwork{}

\makeatletter
\newtheoremstyle{nonumberplainof}%
{\item[\theorem@headerfont\hskip\labelsep ##1\theorem@separator]}%
{\item[\theorem@headerfont\hskip \labelsep ##1\ of\ ##3\theorem@separator]}%
\newtheoremstyle{nonumberplainvon}%
{\item[\theorem@headerfont\hskip\labelsep ##1\theorem@separator]}%
{\item[\theorem@headerfont\hskip \labelsep ##1\ von\ ##3\theorem@separator]}
\makeatother

\theoremstyle{nonumberplainof}
\theoremprework{%
    \renewcommand{\qed}{%
        \hspace*{\fill}\nolinebreak[3]%
        \nopagebreak[3]\hspace*{\fill}{$\square$}}}
\theorempostwork{}

\newtheorem{Proof}{Proof}

\theoremstyle{nonumberplainvon}
\theoremprework{%
    \renewcommand{\qed}{%
        \hspace*{\fill}\nolinebreak[3]%
        \nopagebreak[3]\hspace*{\fill}{$\square$}}}
\theorempostwork{}

\newcommand{\pers}[1]{{\scshape#1}}

\makeatletter
\newcommandx*\newhistpers[3][3]{%
  \ifthenelse{\equal{#3}{}}{%
    \expandafter\newcommand%
    \csname\expandafter\@gobble\string#1\endcsname[1][]{%
      \ifthenelse{\equal{##1}{}}{%
        \pers{#2}%
      }{%
        \pers{#2}##1%
      }%
    }%
  }{%
    \expandafter\newcommand%
    \csname\expandafter\@gobble\string#1\endcsname[1][]{%
      \ifthenelse{\equal{##1}{}}{%
        \pers{#2}%
      }{%
        \pers{#3}##1%
      }%
    }%
  }%
}
\makeatother

\newhistpers{\PBanach}{Banach}
\newhistpers{\PBorel}{Borel}
\newhistpers{\PBregman}{Bregman}
\newhistpers{\PCaratheodory}{Carathéodory}
\newhistpers{\PCauchy}{Cauchy}
\newhistpers{\PDescartes}{Descartes}[Cartes]
\newhistpers{\PDirichlet}{Dirichlet}
\newhistpers{\PDirac}{Dirac}
\newhistpers{\PEuclid}{Euclid}
\newhistpers{\PEuler}{Euler}
\newhistpers{\PFejer}{Fejér}
\newhistpers{\PFenchel}{Fenchel}
\newhistpers{\PFermat}{Fermat}
\newhistpers{\PFredholm}{Fredholm}
\newhistpers{\PFresnel}{Fresnel}
\newhistpers{\PFourier}{Fourier}
\newhistpers{\PFrobenius}{Frobenius}
\newhistpers{\PGateaux}{Gâ\-teaux}
\newhistpers{\PGauss}{Gauß}
\newhistpers{\PGram}{Gram}
\newhistpers{\PHadamard}{Hadamard}
\newhistpers{\PHankel}{Hankel}
\newhistpers{\PHausdorff}{Hausdorff}
\newhistpers{\PHermite}{Hermite}[Hermit]
\newhistpers{\PHesse}{Hesse}
\newhistpers{\PHilbert}{Hilbert}
\newhistpers{\PHoelder}{Hölder}
\newhistpers{\PHurwitz}{Hurwitz}
\newhistpers{\PJacobi}{Jacobi}
\newhistpers{\PJensen}{Jensen}
\newhistpers{\PKronecker}{Kronecker}
\newhistpers{\PLagrange}{Lagrange}
\newhistpers{\PLanczos}{Lanc\-zos}
\newhistpers{\PLaplace}{Laplace}
\newhistpers{\PLebesgue}{Lebesgue}
\newhistpers{\PLegendre}{Legendre}
\newhistpers{\PLidskii}{Lidskii}
\newhistpers{\PLipschitz}{Lipschitz}
\newhistpers{\PMoivre}{Moivre}
\newhistpers{\PMoore}{Moore}
\newhistpers{\PNewton}{Newton}
\newhistpers{\PPaley}{Paley}
\newhistpers{\PParseval}{Parseval}
\newhistpers{\PPenrose}{Penrose}
\newhistpers{\PPoisson}{Poisson}
\newhistpers{\PProny}{Prony}
\newhistpers{\PPythagoras}{Pythagoras}[Pythagor]
\newhistpers{\PRayleigh}{Rayleigh}
\newhistpers{\PRiemann}{Riemann}
\newhistpers{\PRitz}{Ritz}
\newhistpers{\PSard}{Sard}
\newhistpers{\PSchatten}{Schat\-ten}
\newhistpers{\PSchauder}{Schauder}
\newhistpers{\PSchmidt}{Schmidt}
\newhistpers{\PSchwarz}{Schwarz}
\newhistpers{\PSchwartz}{Schwartz}
\newhistpers{\PSobolev}{Sobolev}
\newhistpers{\PTaylor}{Tay\-lor}
\newhistpers{\PTikhonov}{Tik\-ho\-nov}
\newhistpers{\PTitchmarsh}{Titchmarsh}
\newhistpers{\PToeplitz}{Toeplitz}
\newhistpers{\PVandermonde}{Vandermonde}
\newhistpers{\PVieta}{Vieta}
\newhistpers{\PWeierstrass}{Weierstraß}
\newhistpers{\PWiener}{Wiener}
\newhistpers{\PYoung}{Young}

\newcommand{\ie}{i.e.}
\newcommand{\etal}{et~al.}

\newcommand{\Vek}[1]{\bm{#1}}
\newcommand{\Mat}[1]{\bm{#1}}

\newcommand{\Set}[1]{\mathfrak{#1}}
\newcommand{\Op}[1]{\mathop{\kern0pt\mathcal{#1}}\nolimits}
\newcommand{\id}[1]{\mathrm{#1}}

\newcommand{\BR}{\mathbb{R}}

\newcommand{\BN}{\mathbb{N}}

\DeclareMathOperator*{\argmin}{argmin}

\DeclareMathOperator{\diag}{diag}

\DeclareMathOperator{\proj}{proj}
\DeclareMathOperator{\prox}{prox}
\DeclareMathOperator{\ran}{ran}
\DeclareMathOperator{\rank}{rank}

\DeclareMathOperator{\tr}{tr}

\newcommand{\T}{\mathrm{T}}

\newcommand{\iProd}[2]{\left\langle#1,#2\right\rangle}

\newcommand{\iProdn}[2]{\langle#1,#2\rangle}
\newcommand{\iProdb}[2]{\bigl\langle#1,#2\bigr\rangle}

\newcommand{\absn}[1]{\lvert\hspace{1pt}#1\hspace{1pt}\rvert}

\newcommand{\pNormn}[1]{\lVert\hspace{1pt}#1\hspace{1pt}\rVert}

\newlength{\subalignskip}
\newlength{\subalignaboveskip}
\newlength{\subalignbelowskip}
\newlength{\fsmallskip}
\newlength{\fskip}
\newlength{\subintertextskip}
\newlength{\addintertextskip}
\newcommand{\bq}[1]{`#1'}

\newcommand{\tT}{\mathrm{T}}

\AUTHOR{Robert \pers{Beinert}\Inst{1} and Gabriele
  \pers{Steidl}\Inst{1}}
\SHORTAUTHOR{R.~Beinert and G.~Steidl}
\TITLE{Robust PCA via Regularized \pers{reaper} \\[1ex]
with a Matrix-Free Proximal Algorithm}
\SHORTTITLE{Matrix-Free Algorithms for Robust PCA}
\INSTITUTE{\Inst{1}\parbox[t]{0.95\linewidth}{
TU Berlin\\
Straße des 17. Juni 136\\ 
10623 Berlin  
  }
}
\CORRESPONDENCE{R. \pers{Beinert}: \email{beinert@math.tu-berlin.de}\\
  G. \pers{Steidl}: \email{steidl@math.tu-berlin.de}}
\ABSTRACT{
Principal component analysis (PCA) 
is known to be 
sensitive to outliers, so that various robust PCA variants 
were proposed in the literature. A recent model, called \pers{reaper}, 
aims to find the principal components by solving a convex optimization problem.
Usually the number of principal components must be determined in advance and the minimization
is performed over symmetric positive semi-definite matrices 
having the size of the data, although the number of principal components is substantially smaller.
This prohibits its use if the dimension of the data is large which is often the case in image processing.

In this paper, we propose a regularized version of \pers{reaper} which enforces the sparsity
of the number of principal components by penalizing the nuclear norm of the corresponding
orthogonal projector. This has the advantage that only an upper bound on the number of principal components
is required. 
Our second contribution is a matrix-free algorithm to find a minimizer of the regularized \pers{reaper} 
which is also suited for high dimensional data. 
The algorithm couples a primal-dual minimization approach 
with a thick-restarted Lanczos process.
As a side result, we discuss the topic of the bias in robust PCA.
Numerical examples demonstrate the performance of our algorithm.
}
\AMSCLASS{58C05, 62H25, 65K10}

\begin{document}
\thispagestyle{plain}
\TitleHeader

\section{Introduction}\label{sec:introduction}
Principal component analysis (PCA) \cite{Pearson1901} realizes the
dimensionality reduction of data by projecting them onto those affine
subspace which minimizes the sum of the squared Euclidean distances
between the data points and their projections.  Unfortunately, PCA is
very sensitive to outliers, so that various robust approaches were
developed in robust statistics \cite{HR2009,RL1987,Tyler1987} and
nonlinear optimization.  In this paper, we focus on the second one.

One possibility to make PCA robust consists in removing outliers
before computing the principal components which has the serious
drawback that outliers are difficult to identify and other data points
are often falsely labeled as outliers.  Another approach assigns
different weights to data points based on their estimated relevance,
to get a weighted PCA~\cite{KKSZ08} or repeatedly estimate the model
parameters from a random subset of data points until a satisfactory
result indicated by the number of data points within a certain error
threshold is obtained~\cite{fischler1987random}.  In a similar vein,
least trimmed squares PCA models
\cite{podosinnikova2014robust,rousseeuw2005robust} aim to exclude
outliers from the squared error function, but in a deterministic way.
The variational model in \cite{candes2011robust} decomposes the data
matrix into a low rank and a sparse part.  Related approaches such as
\cite{MT2011,XCS2012} separate the low rank component from the column
sparse one using different norms in the variational model.  Another
group of robust PCA replaces the squared $L_2$ norm in the PCA model
by the $L_1$ norm \cite{ke2005robust}.  Unfortunately, this norm is
not rotationally invariant, i.e., when rotating the centered data
points, the minimizing subspace is not rotated in the same way.
Replacing the squared Euclidean norm in the PCA model by just the
Euclidean one, leads to a non-convex robust PCA model with
minimization over the Stiefel or Grassmannian manifold, see,
e.g.~\cite{ding2006r,LM2018,MZL2019,NNSS2020}.  Instead of the
previous model which minimizes over the sparse number of directions
spanning the low dimensional subspace, it is also possible to minimize
over the orthogonal projectors onto the desired subspace.  This has
the advantage that the minimization can be performed over symmetric
positive semi-definite matrices, e.g. using methods from semi-definite
programming, and the disadvantage that the dimension of the projectors
is as large as the data now. This prohibits this approach for many
applications in particular in image processing.  The projector PCA
model is still non-convex and a convex relaxation, called
\pers{reaper}, was recently proposed by \pers{Lerman} et
el. \cite{LCTZ15}.

In this paper, we build up on the advantages of the convex \pers{reaper} model, but modify it in two 
important directions:
(i) by penalizing the nuclear norm of the approximated projectors, our model does only require
an upper bound on the dimension of the desired subspace.
 Having the same effect as the sparsity promotion of the
1-norm, the nuclear norm -- the 1-norm of the eigenvalues -- promotes low-rank matrices or, equivalently, sparse
eigenvalue decompositions;
(ii) 
by combining primal-dual minimization techniques with a thick-restarted Lanczos process,
we are able to handle high dimensional data.
We call our new model r\pers{reaper}. We provide all computation steps leading to a provable convergent algorithm
and give a performance analysis following the lines of \cite{LCTZ15}.
The choice of the offset in robust PCA is an interesting problem which is not fully discussed in the literature so far.
Usually, the geometric median is used. We do not provide a full solution of this issue, 
but show that under some assumptions
the affine hyperplane in $\mathbb R^d$
having the smallest Euclidean distance to $n > d$ given data points goes through $d+1$  of these points.
We underline our theoretical findings by numerical examples.

The outline of this paper is as follows:
preliminaries from linear algebra and convex analysis are given in 
in Section \ref{sec:notation}. In Section \ref{sec:model},
we introduce our regularized \pers{reaper} model.
The basic primal-dual algorithm for its minimization is 
discussed in Section \ref{sec:primal-dual-method}.
The algorithm is formulated with respect to the full 
projection matrix. The matrix-free version of the algorithm
is given in Section \ref{sec:matrix-free-real}.
It is based on the thick-restarted Lanczos algorithm and
is suited for high-dimensional data.
In Section \ref{sec:error-analysis}, we 
examine the performance analysis of r\pers{reaper} along the lines of  \cite{LCTZ15}.
Some results on the offset in robust PCA are proved in  Section \ref{sec:incorporating-offset}.
The very good performance of r\pers{reaper} in particular for high dimensional data is 
demonstrated in Section \ref{sec:numerics}.
Section \ref{sec:conclusion} finishes the paper with conclusions and
directions of future research.

\section{Notation and preliminaries}\label{sec:notation}
Throughout this paper we will use the following notation and basic facts from linear algebra and convex analysis
 which can be found in detail in various monographs and overview papers as
\cite{Bec17,BSS16,CP16,GV13,Ro1970}.

\paragraph{Linear algebra.}
By $\| \cdot\|_2$ we denote the Euclidean vector norm and 
by  $\| \cdot\|_1$ the norm which sums up the absolute vector components.
Recall that for any $\Vek x \in \mathbb R^n$,
\begin{equation}\label{mittel}
\tfrac{1}{\sqrt{n}} \| \Vek x\|_1 \le \| \Vek x\|_2 \le \| \Vek x\|_2.
\end{equation}
Let $\mathbf{1}_n$ resp. $\mathbf{0}_n$ be the vectors having $n$ entries $1$, resp., $0$. 
Analogously, we write $\Mat 1_{n,d}$ and
  $\Mat 0_{n,d}$ for the all-one and all-zero matrix in
  $\BR^{n,d}$. Further, $\Mat I_n$ is the
$n \times n$ identity matrix.  Let $\tr \Mat A$ denote the
\emph{trace} of the quadratic matrix $\Mat A \in \mathbb R^{n,n}$,
i.e., the sum of its eigenvalues.  On $\mathbb R^{n,d}$
the Hilbert--Schmidt inner product is defined by
$$
\langle \Mat X,\Mat Y \rangle \coloneqq \tr(\Mat X^\tT \Mat Y) = \tr(\Mat Y \Mat X^\tT), \qquad \Mat X,\Mat Y \in \mathbb R^{n,d},
$$
and the corresponding so-called \emph{Frobenius norm} by
$\|\Mat X\|_F^2 = \langle \Mat X , \Mat X \rangle$.

Let $\mathcal S(n) \subset \mathbb R^{n,n}$ denote the linear subspace
of symmetric matrices.  For two symmetric matrices $\Mat A,\Mat B \in \mathcal S(n)$, we write
$\Mat A \preceq \Mat B$ if $\Mat B - \Mat A$ is positive
semi-definite.  Every $\Mat A \in \mathcal S(n)$ has a spectral
decomposition
$$
\Mat A = \Mat U \diag(\Vek \lambda_{\Mat A}) \, \Mat U^\tT,
$$
where $\Vek \lambda_{\Mat A} \in \mathbb R^n$ denotes the vector containing the eigenvalues of $\Mat A$ in descending order
$\lambda_1 \ge \cdots \ge \lambda_n$ and $\Mat U$ is the orthogonal matrix
having the corresponding orthogonal eigenvectors as columns. 
The nuclear norm (trace norm) of $\Mat A \in \mathcal S(n)$ is given by
$$
\|\Mat A\|_{\tr} \coloneqq \sum_{j=1}^n |\lambda_j|.
$$
The trace and Frobenius norm correspond to the Schatten 1-norm and
2-norm respectively, where the Schatten $p$-norm with
$1 \le p \le \infty$ of a symmetric matrix $\Mat A$ is defined by
$\pNormn{\Mat A}_{S_p} \coloneqq \pNormn{\Vek \lambda_{\Mat A}}_p$.
Recall that $\Mat \Pi \in \mathbb R^{n,n}$ is an orthogonal projector if
$\Mat \Pi \in \mathcal S(n)$ and $\Mat \Pi^2 = \Mat \Pi$.  This is
equivalent to the statement that $\Mat \Pi \in \mathcal S(n)$ and has
only eigenvalues in $\{0,1\}$.  The \emph{nuclear norm} is the unique
norm such that
\begin{equation*}
\mathrm{rank} (\Mat \Pi) = \|\Mat \Pi\|_{\tr}
\end{equation*}
for every orthogonal projector $\Mat \Pi$.

For a given norm $\| \cdot \|$ on $\mathbb R^n$, the \emph{dual norm} is defined by
$$
\|\Vek x\|_* \coloneqq \max_{\|\Vek y\| \le 1}  \langle \Vek x,\Vek y \rangle.
$$
In particular, for a matrix $\Mat X = (\Vek x_1|\ldots|\Vek x_N) \in \mathbb R^{n,N}$ we will be interested in the norm
$$
\pNormn{\Mat X}_{2,1} \coloneqq \sum_{k=1}^N \pNormn{\Vek x_k}_2
$$
which can be considered as norm on $\mathbb R^{nN}$ by arranging the columns of the matrix into a vector.
Its dual norm is given by 
\begin{equation*} 
\pNormn{\Mat X}_{2,1,*} = \pNormn{\Mat X}_{2,\infty} \coloneqq \max_{k=1,\ldots,N} \pNormn{\Vek x_k}_2.
\end{equation*}

\paragraph{Convex analysis.}
Let $\Gamma_0(\mathbb R^n)$ denote the space of proper, lower semi-continuous, convex functions
mapping from $\mathbb R^n$ into the extended real numbers $(-\infty,\infty]$.
The \emph{indicator function} $\iota_{\mathcal C}$ of $\mathcal C
\subseteq \mathbb R^n$ is defined by
$$
\iota_{\mathcal C}(\Vek x) = 
\left\{
\begin{array}{ll}
0&\mathrm{if} \; \Vek x \in \mathcal C,\\
+ \infty&\mathrm{otherwise}.
\end{array}
\right.
$$
We have $\iota_{\mathcal C} \in \Gamma_0(\mathbb R^n)$ if and only if
$\mathcal C$ is non-empty, convex and closed.

For $f \in \Gamma_0(\mathbb R^n)$, the \emph{proximal mapping} is defined by
\begin{equation*}
  \prox_{f}(\Vek x) \coloneqq \argmin_{\Vek y \in \BR^n} \left\{f(\Vek y) +
  \tfrac12 \pNormn{\Vek x - \Vek y}_2^2 \right\}.
\end{equation*}
Indeed, the minimizer exists and is unique \cite[Thm~31.5]{Ro1970}.
If $\mathcal C \subset \mathbb R^n$ is a nonempty, closed, convex set, 
then the proximal mapping of a multiple of $\iota_{\mathcal C}$ is just
the orthogonal projection onto $\mathcal C$, i.e.,
$$
\prox_{\sigma \iota_{\mathcal C}}(\Vek x) = \proj_{\mathcal C} (\Vek x), \qquad \sigma >0.
$$
In particular, the orthogonal projection onto the halfspace
$\mathcal H(\Vek a,\beta) \coloneqq \{\Vek x \in \mathbb R^n: \langle \Vek
a , \Vek x \rangle \le \beta\}$ with $\Vek a \in \BR^{n}$ and $\beta \in \BR$
can be computed by
\begin{equation*}
\proj_{\mathcal H(\Vek a,\beta)} (\Vek x) = \Vek x - \frac{(\langle \Vek
  a, \Vek x \rangle - \beta)_+}{\|\Vek a\|_2^2} \, \Vek a,
\end{equation*}
where $(y )_+ \coloneqq \max\{0,y\}$. Further, the
orthogonal projection onto the hypercube $Q := [0,1]^n$ is given by 
\begin{equation}\label{proj_cube}
\proj_{\mathcal Q} (\Vek x) = \left(  \max \left\{ \min\{ x_j,1 \},0 \right\} \right)_{j=1}^n.
\end{equation}

The \emph{Fenchel dual} of $f \in \Gamma_0(\mathbb R^n)$ is the function $f^* \in \Gamma_0(\mathbb R^n)$ 
defined by
$$
f^*(\Vek p) \coloneqq \max_{\Vek x \in \mathbb R^n}\langle \Vek p,\Vek x \rangle - f(\Vek x).
$$ 
The dual function of a norm is just the indicator function of the unit ball with respect to its dual norm.
In particular, we have for $\| \cdot\|_{2,1}: \mathbb R^{n,N} \rightarrow \mathbb R$ that
\begin{equation} \label{spec_fenchel} 
\|\Mat X \|_{2,1}^* = \iota_{\mathcal B_{2,\infty}}(\Mat X),
\end{equation}
where 
$
\mathcal B_{2,\infty} \coloneqq \{\Mat X \in \mathbb R^{n,N}: \| \Vek x_k \|_2 \le 1 \; \mathrm{for \; all} \; k=1,\ldots,N\}
$.

\section{Regularized \pers{reaper}} \label{sec:model}
Given $N$ data points $\Vek x_1,\ldots,\Vek x_N \in \BR^n$, the classical
PCA finds a $d$-dimensional affine subspace
$\{\Mat A \, \Vek t +  \Vek b: \Vek t \in \BR^d\}$, $1 \le d \ll n$,
by minimizing
\begin{equation} \label{PCA_1}
         \sum_{k=1}^N \min_{t \in \BR^d} \|\Mat A\, \Vek t + \Vek b - \Vek x_k \|_2^2 
        = 
         \sum_{k=1}^N \|(\Mat A \Mat A^\tT - \Mat I_n )(\Vek b-\Vek x_k)\|_2^2 
                                \quad \mathrm{subject \; to} \quad \Mat A^\tT \Mat A = \Mat I_d
\end{equation}
over $\Vek b \in \BR^n$ and $\Mat A \in \BR^{n,d}$. %
It is not hard to check that the affine subspace 
goes through the \emph{offset} (bias) 
\begin{equation}\label{mean}
\bar {\Vek b} \coloneqq \tfrac{1}{N}(\Vek x_1 + \ldots + \Vek x_N).
\end{equation}
Therefore, we can reduce our attention to data points
$\Vek x_k - \bar{\Vek b}$, $k=1,\ldots,N$, which we denote by
$\Vek x_k$ again, and minimize over the linear $d$-dimensional
subspaces through the origin, i.e.,
\begin{equation*}
  \min_{\Mat A \in \mathbb R^{n,d}} \sum_{k=1}^N \|(\Mat A \Mat A^\tT - \Mat I_n ) \Vek x_k \|_2^2 
  \quad \mathrm{subject \; to} \quad \Mat A^\tT \Mat A = \Mat I_d,
\end{equation*}
where $\Mat X \coloneqq (\Vek x_1| \ldots |\Vek x_N) \in \BR^{n,N}$.

Unfortunately, the solution of this minimization problem is sensitive
to outliers. Therefore several robust PCA variants were proposed in
the literature.  A straightforward approach consists in just skipping
the square in the Euclidean norm leading
to
\begin{equation} \label{PCA_rob}
\min_{\Mat A \in \mathbb R^{n,d}} \sum_{k=1}^N \|(\Mat A \Mat A^\tT -
\Mat I_n ) \Vek x_k \|_2
= \| \Mat A \Mat A^\tT \Mat X - \Mat X\|_{2,1} \quad \mathrm{subject \; to} \quad \Mat A^\tT \Mat A = \Mat I_d.
\end{equation}
This is a nonconvex model which requires the minimization over
matrices $\Mat A$ in the so-called Stiefel manifold, 
$$
\mathrm{St}(n,d) \coloneqq \{\Mat A \in \BR^{n,d}: \Mat A^\tT \Mat A = \Mat I_d \}.
$$ 
see \cite{ding2006r,LM2018,NNSS19b,NNSS2020}.

Another approach is based on the observation that
$\Mat \Pi \coloneqq \Mat A \Mat A^\tT$ is the orthogonal projector onto
the linear subspace spanned by the columns of $\Mat A$.  Since the
linear subspace is $d$-dimensional, exactly $d$ eigenvalues of
$\Mat \Pi$ have to be one.  Thus, problem \eqref{PCA_rob} can be
reformulated as
 \begin{equation} \label{PCA_2}
        \min_{\Pi \in \mathcal S(n) } \|\Mat \Pi \Mat X - \Mat X \|_{2,1} 
                                \quad \mathrm{subject \; to} \quad
                                \Vek \lambda_{\Mat \Pi} \in \{0,1\}^n,
                                \; \tr(\Mat \Pi) = d. 
\end{equation}
Having computed $\Mat \Pi$, we can determine $\Mat A$ by spectral decomposition.
Unfortunately, \eqref{PCA_2} is still a nonconvex model which is moreover  NP hard to solve. Therefore
Lerman et al. \cite{LCTZ15} suggested to replace it by a convex relaxation, called \pers{reaper},
\begin{equation*} 
        \min_{P \in \mathcal S(n) } \|\Mat P \Mat X - \Mat X \|_{2,1} \quad \mathrm{subject \; to} 
        \quad \Mat 0_{n,n} \preceq  \Mat P \preceq \Mat I_n, \; \tr(\Mat P) = d.
\end{equation*}
In order to deal with the
non-differentiability of the objective function, Lerman \etal\ \cite{LCTZ15}
iteratively solve a series of positive semi-definite programs.
In contrast to models minimizing directly over
$\Mat A \in \mathbb R^{n,d}$, algorithms for minimizing \pers{reaper}
or r\pers{reaper} seem to require the handling of a large matrix
$\Mat P \in \mathcal S(n)$ or, more precisely, the
  handling of its spectral decomposition which makes
the method not practicable for high-dimensional data. 

The above model requires the exact knowledge of the dimension $d$ of the linear subspace the data will be reduced to.
In this paper, we suggest to replace the strict trace constraint by a relaxed variant $\tr(\Mat \Pi) \le d$ 
and to add the nuclear norm of $\Mat \Pi$ as a regularizer
which enforces the sparsity of the rank of $\Mat \Pi$:
\begin{equation} \label{PCA_3_1}
\min_{\Pi \in \mathcal S(n) } \|\Mat \Pi \Mat X - \Mat X \|_{2,1}
        + \alpha \|\Mat \Pi\|_{\tr}  
        \quad \mathrm{subject \; to} \quad \Vek \lambda_{\Mat \Pi} \in
        \{0,1\}^n, \; \tr(\Mat \Pi) \le d. 
\end{equation}
Here $\alpha >0$ is an appropriately fixed regularization parameter.

Since \eqref{PCA_3_1} is again hard so solve, we use a relaxation for the eigenvalues and
call the new model regularized \pers{reaper} (r\pers{reaper}):
\begin{equation} \label{PCA_3}
        \min_{P \in \mathcal S(n) } \|\Mat P \Mat X - \Mat X \|_{2,1}
        + \alpha \|\Mat P\|_{\tr}  
        \quad \mathrm{subject \; to} \quad \Mat 0_{n,n} \preceq  \Mat P \preceq \Mat I_n, \; \tr(\Mat P) \le d.
\end{equation}
Finally, we project the solution of r\pers{reaper} to the set of orthoprojectors with rank not larger than $d$:
\begin{equation*}
\mathcal O_d \coloneqq \{ \Mat \Pi \in \mathcal S(n): \Vek \lambda_{\Vek \Pi} \in \mathcal E_d\},
\end{equation*}
where
\begin{equation*}
\mathcal E_d \coloneqq \{ \lambda \in \BR^n: \Vek \lambda \in \{0,1\}^n, \, \langle \Vek \lambda, \Vek 1_n \rangle \le d\}.
\end{equation*}

In the following we will present a primal-dual
approach to solve \eqref{PCA_3} which uses only the sparse
spectral decomposition of $\Mat P$, but not the matrix itself within
the computation steps.

\section{Primal-dual algorithm} \label{sec:primal-dual-method}
r\pers{reaper} 
is  a convex optimization problem; so
we may choose from various convex solvers.  Since both -- data
fidelity and nuclear norm -- are non-differentiable, we apply
the primal-dual method of Chambolle and Pock
\cite{CP16}.  
For this purpose, we define the forward operator
$$\Op X \colon \mathcal S(n) \to \BR^{n, N} : \Mat P \mapsto \Mat P \Mat X$$
and rearrange \eqref{PCA_3} as
\begin{equation}  \label{eq:reaper-matrix-form}
  \min_{\Mat P \in \mathcal S(n)} 
    \pNormn{\Op X (\Mat P) - \Mat X}_{2,1}
  + \alpha \Op R (\Mat P) ,
 \end{equation}
where the regularizer
$\Op R \colon \mathcal S(n) \to [0,+\infty]$ is defined by
\begin{equation}  \label{eq:def-reg}
  \Op R(\Mat P)
  \coloneqq 
  \pNormn{\Mat P}_{\tr} + \iota_{\mathcal C}(\Mat P),
  \qquad
  \mathcal C
  \coloneqq
  \{
  \Mat P \in \mathcal S(n)
  :
  \Mat 0_{n,n} \preceq \Mat P \preceq \Mat I_n,
  \tr( \Mat P) \le d
  \}.
\end{equation}
Since $\mathcal C$ is compact and convex, and since the
norms $\pNormn{\cdot}_{2,1}$ and $\pNormn{\cdot}_{\tr}$
are continuous, r\pers{reaper} has a global
minimizer.  This minimizer is in general not unique.  Concerning the
adjoint operator $\Op X^*: \BR^{n, N} \to \mathcal S(n)$, we
observe
\begin{align*}
  \iProdn{\Op X(\Mat P)}{\Mat Y}
  &=
    \tfrac12 \bigl( \iProdn{\Mat P \Mat X}{\Mat Y} + \iProdn{\Mat P^\T
    \Mat X}{\Mat Y} \bigr)\\
  & =
    \tfrac12 \bigl( \tr( \Mat Y^\T \Mat P \Mat X) + \tr(\Mat X^\T \Mat P
    \Mat Y) \bigr)
    =
    \tfrac12 \iProd{\Mat P}{\tfrac12(\Mat X \Mat Y^\T + \Mat Y \Mat X^\T) }.
\end{align*}
for all $\Mat P \in \mathcal S(n)$ and
  $\Mat Y \in \BR^{n,N}$, where we exploit the symmetry of $\Mat P$ by
  $\Mat P = \nicefrac12(\Mat P + \Mat P^\T)$.  Thus,
the adjoint is just
\begin{equation*}
  \Op X^*(\Mat Y) = \tfrac12(\Mat X \Mat Y^\T + \Mat Y \Mat X^\T).
\end{equation*}

The operator norm of $\Op X$ is given by the spectral norm of
$\Mat X \in \BR^{n,N}$, \ie\ 
$$\| \Op X\| = \|\Mat X\|_2.$$
In more detail, for $\Mat P = (\Vek p_1|\ldots|\Vek
  p_n) \in \mathcal S(n)$, we obtain
\begin{align*}
  \| \Op X\| = \max_{\Mat P \in \mathcal S(n) \atop \|\Mat P\|_F \le
  1} \|\Mat P \Mat  X\|_F  
  = 
  \max_{\Mat P \in \mathcal S(n) \atop \|\Mat P\|_F \le 1}
  \biggl(\sum_{j=k}^n \|\Mat X^\tT \Vek p_k\|_2^2\biggr)^\frac12 
  \le 
  \max_{\Mat P \in \mathcal S(n) \atop \|\Mat P\|_F \le 1}
  \biggl(\|\Mat X\|_2^2 \sum_{k=1}^n \|\Vek p_k\|_2^2\biggr)^\frac12 
  \le
  \|\Mat X\|_2.
\end{align*}             
Here the inequality becomes sharp for
$\Mat P = \Mat U \diag ( (1,0,\ldots,0)^\tT) \, \Mat
U^\tT$, where $\Mat U$ arises from the singular value decomposition
$\Mat X = \Mat U \diag( \Vek \sigma_{\Mat X}) \, \Mat V^\tT$ with descending ordered
singular values $\sigma_1 \ge \cdots \ge \sigma_{\min\{n,N\}}$.

Next, we apply the primal-dual method of Chambolle and
Pock \cite{CP16} with extrapolation of the primal variable to compute
the minimizer of r\pers{reaper}
\eqref{eq:reaper-matrix-form}, which leads us to the
  following numerical method.

\begin{Algorithm}[Primal-Dual Algorithm]   \label{alg:proj-trunc-cube_0}
  {\scshape Input:} $\Mat X \in \BR^{n,N}$ $d \in \BN$, and $\sigma,\tau > 0$ with $\sigma\tau < 1/\|\Mat X\|_2^2$, and
        $\theta \in (0,1]$.\\
  {\scshape Intialization:} ${\Mat P}^{(0)} = \bar{\Mat P}^{(0)} = \Mat 0_{n,n}$ ,
        $\Mat Y^{(0)} \coloneqq \Mat 0_{n, N}$.
        \\
        {\scshape Interation:}
        \begin{align*}
         \Mat Y^{(r+1)}
  &\coloneqq \prox_{\sigma \pNormn{\cdot \, - \Mat X}_{2,1}^*} \Bigl(
    \Vek Y^{(r)} + \sigma \Op 
    X \bigl(\bar{\Mat P}^{(r)} \bigr) \Bigr),
  \\[\fskip]
  \Mat P^{(r+1)}
  &\coloneqq \prox_{\tau \alpha \Op R} \Bigl( \Vek P^{(r)} - \tau \Op
    X^* \bigl(\Vek Y^{(r+1)} \bigr) \Bigr),
  \\[\fskip]
  \bar{\Mat P}^{(r+1)}
  &\coloneqq (1 + \theta) \, \Mat P^{(r+1)} - \theta \, \Mat P^{(r)}
\end{align*} 
\end{Algorithm}

More generally, Chambolle and Pock \cite{CP16} have
  proven that the sequence
$\{\Mat P^{(r)} \}_{r\in\BN}$ converges to a minimizer $\hat {\Mat P}$ of
\eqref{eq:reaper-matrix-form} and the sequence $\{\Mat Y^{(r)} \}_{r\in\BN}$
to a minimizer of the dual problem
$$
\min_{\Mat Y \in \mathbb R^{n,N}} \| \cdot - X\|_{2,1}^*(\Mat Y) + (\alpha \mathcal R)^*(- \Op X^*(\Mat Y) )
$$
if the Lagrangian 
$$
L(\Mat P,\Mat Y) \coloneqq - \| \cdot - X\|_{2,1}^*(\Mat Y) + \alpha \mathcal R(\Mat P) + \langle \Op X(\Mat P), \Mat Y \rangle
$$
has a saddle-point which is, however, clear
for r\pers{reaper}.

The algorithm requires the computation of the proximal  mapping of the dual data fidelity
and of the regularizer which we consider next.

\begin{Proposition}[Proximal mapping of the dual data fidelity]  \label{thm:prox-spec}
For $\Mat x \in \mathbb R^{n,N}$ and $\sigma > 0$, we have
$$
 \prox_{\sigma \pNormn{\cdot \, - \Mat X}_{2,1}^*} = \proj_{\Set B_{2,\infty}} ( \cdot - \sigma \Mat X).
$$
\end{Proposition}

\begin{Proof}
 Using \eqref{spec_fenchel} and, since
  $(f(\cdot - x_0))^* = f^* + \langle \cdot, x_0 \rangle$, we obtain
  \begin{align*}
    \prox_{\sigma \pNormn{\cdot \, - \Mat X}_{2,1}^*} (\Mat Y)
    &=
      \argmin_{\Mat Z \in \mathbb R^{n,N}} \bigl\{\tfrac12 \|\Mat Z-\Mat Y\|_F^2 
      +  \iota_{\Set B_{2,\infty}}(\Mat Z) + \sigma \langle \Mat
      Z,\Mat X \rangle  \bigr\}\\ 
    &=
      \argmin_{\Mat Z \in \mathbb R^{n,N}} \bigl\{ \tfrac12 \|\Mat
      Z-(\Mat Y-\sigma \Mat X)\|_F^2 +  \iota_{\Set B_{2,\infty}}(\Mat
      Z)\bigr\}\\ 
    &= \proj_{\Set B_{2,\infty}} ( Y - \sigma \Mat X).
      \tag*{\qed}
  \end{align*}
\end{Proof}

For the maximal dimension $d$ of the target subspace,
  we henceforth use the half-space
$$
\mathcal{H} \coloneqq \mathcal{H}(\Vek 1_n,d) = \{\Vek x \in \mathbb R^n: \langle \Vek x,\Vek 1-n \rangle \le d \}.
$$
in order to bound the trace of the primal iteration variable
  $\Mat P^{(r)}$.  Then the proximal mapping of the regularizer is given in the following proposition.

\begin{Proposition}[Proximal mapping of the regularizer]  \label{thm:prox-spec}
For  $\Mat P \in \mathcal S(n)$ 
with spectral decomposition 
$\Mat P = \Mat U \diag(\Vek \lambda_{\Mat P}) \, \Mat U^\tT$ and $\mathcal R$ in \eqref{eq:def-reg}
it holds
\begin{equation*}
    \prox_{\tau \alpha \Op R}(\Mat P)
    =
    \Mat U \diag(\proj_{\mathcal{Q}\cap \mathcal{H}}(\Vek \lambda_{\Mat P} - \tau \alpha \Vek 1_n)) \, \Mat U^\tT. \tag*{\qed}
  \end{equation*}
\end{Proposition}

\begin{Proof}
A symmetric matrix $\Mat P$ is in $\mathcal C$ if and only if
$\Vek \lambda_{\Mat P} \in \mathcal{Q} \cap \mathcal{H}$.
Hence the regularizer can be written as
\begin{equation*}
  \Op R (\Mat P)  
  = \langle \Vek \lambda_{\Mat P}, \Vek 1_n \rangle
    + \iota_{\mathcal{Q} \cap {\mathcal H}}(\Vek \lambda_{\Mat P}).
\end{equation*}
and
\begin{equation} \label{aa}
\prox_{\tau \alpha \Op R}(\Mat P) 
= \argmin_{\Mat S \in \mathcal S(n)} 
\bigl\{ \tfrac12 \|\Mat S - \Mat P\|_F^2 + \tau \alpha \langle \Vek \lambda_{\Mat S}, \Vek 1_n \rangle
+ \iota_{\mathcal{Q}  \cap {\mathcal H}}(\Vek \lambda_{\Mat S}) \bigr\}.
\end{equation}
By the theorem of Hoffmann and Wielandt \cite[Theorem 6.3.5]{HJ1991}, we
know that
$$
\|\Mat S - \Mat P\|_F^2 \ge \|\Vek \lambda_{\Mat S} - \Vek \lambda_{\Mat P} \|_2^2
$$
with equality if and only if $\Mat S$ possesses the same eigenspaces
as $\Mat P$. 
Therefore, the minimizer in \eqref{aa} has to be of the form
$\Mat S = \Mat U \diag(\Vek \lambda_{\Mat S}) \, \Mat
  U^\tT$, where the columns of $\Mat U$ are the eigenvectors of $\Mat P$.  
Incorporating this observation in (\ref{aa}),
we determine the eigenvalues $\Vek \lambda_{\Mat S}$ by solving the
minimization problem
\begin{align*}
  \Vek \lambda_{\Mat S}
  &= \argmin_{\Vek \lambda_{\Mat S} \in \mathbb R^n}
    \Bigl\{ \tfrac12 \|\Vek \lambda_{\Mat S} - \Vek \lambda_{\Mat P}\|_2^2 + \tau \alpha \langle \Vek \lambda_{\Mat S}, \Vek 1_n \rangle
    + \iota_{\mathcal{Q} \cap {\mathcal H}}(\Vek \lambda_{\Mat S}) \Bigr\}
  \\
  &= \argmin_{\Vek \lambda_{\Mat S}} \Bigl\{ \tfrac12 \|\Vek \lambda_{\Mat S} + \tau \alpha \Vek 1_n - \Vek \lambda_{\Mat P}\|_2^2 
    + \iota_{\mathcal{Q} \cap {\mathcal H}}(\Vek \lambda_{\Mat S})
    \Bigr\}
  \\
  &= 
    \proj_{\mathcal{Q}\cap{\mathcal H}}(\Vek \lambda_{\Mat P} - \tau
    \alpha \Vek 1_n)).
    \tag*{\qed}
\end{align*}
\end{Proof}

Alternatively to the proof we could argue with the so-called spectral function related to $\mathcal R$ which is invariant under permutations, 
see, e.g. \cite{Bec17}.

By Proposition \ref{thm:prox-spec} the proximal mapping of the regularizer
requires the projection onto the truncated hypercube. The following proposition can be found in
\cite[Ex~6.32]{Bec17}.

\begin{Proposition}[projection onto the truncated hypercube]  \label{thm:proj_hyperqube}
For any $\Vek \lambda \in \BR^n$ and any $d \in (0,n]$, the projection
to the truncated hypercube is given by 
$$
\proj_{\mathcal{Q}\cap{\mathcal H}}(\Vek \lambda) =
\left\{
\begin{array}{ll}
\proj_{\mathcal{Q} } (\Vek \lambda)                  &\; \mathrm{if} \; \langle \proj_{\mathcal{Q} } (\Vek \lambda), \Vek 1_n \rangle \le d,\\
\proj_{\mathcal{Q} } (\Vek \lambda - \hat t \Vek 1_n)&\; \mathrm{otherwise} ,
\end{array}
\right.
$$
where $\hat t$ is the positive root of the function
\begin{equation}\label{phi}
\varphi(t) \coloneqq \langle \proj_{\mathcal{Q} }( \Vek \lambda - t \Vek 1_n ), \Vek 1_n \rangle - d.
\end{equation}
\end{Proposition}

Due to the projection to the hypercube, see \eqref{proj_cube}, 
only the positive components of $\Vek \lambda$
influence its projection onto $\mathcal{Q}\cap{\mathcal H}$.  More
precisely, we have
$$\proj_{\mathcal{Q}\cap{\mathcal H}}(\Vek \lambda) = \proj_{\mathcal{Q}\cap{\mathcal H}}(\Vek \lambda)_+ ,$$
where the function $(\cdot)_+$ is employed componentwise.

To formulate a projection algorithm, in particular, to compute the
zero of $\varphi$, we study the
properties of $\varphi$.

\begin{Lemma}[Properties of $\varphi$] \label{lem:varphi} For fixed
  $\Vek \lambda \in \mathbb R^n$ with
    $\iProdn{\proj_{\mathcal Q}(\Vek \lambda)}{\Vek 1_n} > d$, the
  function $\varphi:[0,\infty) \rightarrow \mathbb R$ defined in
  \eqref{phi} has the following properties:
  \begin{itemize}
  \item[i)] $\varphi$ is Lipschitz continuous.
  \item[ii)] There exists $M \in \mathbb N$,
    $M < 2n$ and
    $0 = s_0 < s_1 < s_2 < \ldots < s_M < s_{M+1}$ such that
    $\varphi(t) > 0$ for $t \in [0,s_{M}]$ and $\varphi(t) \le 0$ for
    $t \ge s_{M+1}$.  Further, we have piecewise linearity
    $$
    \varphi(t) = \varphi(s_l) - k_l (t-s_l), \quad t \in [s_l,s_{l+1}), \; l = 0,\ldots, M,
    $$
    where
    $$
    k_l \coloneqq |\{j \in \{1,\ldots,n\}: (\Vek \lambda - s_l \Vek
    1_n)_j \in (0,1] \}|. 
    $$
    In particular, the function $\varphi$ is monotone decreasing.
  \item[iii)] The positive zero $\hat t$ of $\varphi$ is given by
    $$
    \hat t = s_M + \tfrac{1}{k_M} \varphi(s_M). 
    $$
  \end{itemize}
\end{Lemma}

\begin{Proof}
i) Using the definition of $\varphi$,
the Cauchy--Schwarz inequality, and
the nonexpansiveness of the projection, we get 
\begin{align*}
  |\varphi(t) - \varphi(s)| 
  &= |\langle \proj_{\mathcal{Q} }(\Vek \lambda - t \Vek 1_n ), \Vek 1_n \rangle - \langle \proj_{\mathcal{Q} }(\Vek \lambda - s\Vek 1_n ), \Vek 1_n \rangle|\\
  &\le 
    \sqrt{n} \, \|\proj_{\mathcal{Q} }(\Vek \lambda - t \Vek 1_n )  -
    \proj_{\mathcal{Q} }(\Vek \lambda - s\Vek 1_n )\|_2\\ 
  &\le
    \sqrt{n} \, \|(s-t) \Vek 1_n\|_2 = n \, |s-t|.
\end{align*}
ii) By definition of $\varphi$ and by the assumption $\iProdn{\proj_{\mathcal Q}(\Vek \lambda)}{\Vek 1_n} > d$,
we have $\varphi(0) > 0$.
Starting with $s_0 = 0$, we construct $s_l$ with $l=1,\ldots,M$ iteratively as follows: 
given $s_l$ with $\varphi(s_l) >0$, we set $\Vek \mu \coloneqq \Vek \lambda - s_l \Vek 1_n$ and choose
$$
s_{l+1} \coloneqq s_{l} + h_l, \qquad h_l := \min\{s_{\mathrm{leave}},s_{\mathrm{enter}} \},
$$
where
$$
s_{\mathrm{leave}} \coloneqq \min_j \{\mu_j: \mu_j \in (0,1]\}, \quad
s_{\mathrm{enter}} \coloneqq \min_j \{\mu_j-1: \mu_j >1\}.
$$
Here we use the convention $\min \emptyset = \infty$.
Note that at least one of the above sets in the
  definition of $s_{\mathrm{leave}}$ and $s_{\mathrm{enter}}$ is
  non-empty since otherwise all components of $\Vek \mu$
have to be non-positive implying $\proj_{\mathcal{Q} }(\Vek \mu) = \Vek 0_n$ and thus
$\varphi(s_l) = -d$, a contradiction.

Considering the projection to the hypercube
  $\mathcal Q$ in (\ref{proj_cube}), we see that the index set
$\{j \in \{1,\ldots,n\}: (\Vek \lambda - t \Vek 1_n)_j \in (0,1] \}$
does not change for $t \in [s_l, s_{l+1})$ and that a
change appears exactly in
$s_{l+1}$, where at least one component enters or leaves the interval
$(0,1]$.  Hence we have
$$
\varphi(t) = \varphi(s_l) - k_l (t-s_l), \quad t \in [s_l,s_{l+1}).
$$
Let $s_{M+1}$ be the first value in this procedure, where $\varphi(s_{M+1}) \le 0$.
Since each component in $\Vek \lambda - t \Vek 1_n$ can at most one times enter or leave the interval $(0,1]$,
we know that $M < 2n$.
Further, we have $k_M >0$ since our piecewise linear
  function cannot pass zero in the interval $[s_M, s_{M+1}]$ otherwise.

iii) Now the zero $\hat t$ of $\varphi$ in the
interval $[s_{M},s_{M+1}]$ can be computed by solving
\begin{equation*}
\varphi(\hat t) = \varphi( s_M) - k_M (\hat t-s_M) = 0,
\end{equation*} which results in 
$\hat t = s_M + \tfrac{1}{k_M} \varphi(s_M)$ 
and finishes the proof.  \qed
\end{Proof}

Following Proposition \ref{thm:proj_hyperqube} and the previous proof, 
we obtain the following algorithm for the projection onto $\mathcal{Q} \cap \mathcal{H}$.

\begin{Algorithm}[Projection onto truncated hypercube]  \label{alg:proj-trunc-cube}
  {\scshape Input:} $\Vek \lambda \in \BR^n$, $d \in \mathbb N$.
  \begin{enumerate}[(i), nosep]
        \item
  Compute $\Vek \mu \coloneqq \proj_{\mathcal Q} (\Vek \lambda)$ by \eqref{proj_cube}. \\
  If $\langle \Vek \mu, \Vek 1_n \rangle \le d$, then return $\hat{\Vek \lambda} = \Vek \mu$;
        \\
        otherwise set $s \coloneqq  0$, $\varphi \coloneqq + \infty$ and $\Vek \mu = \Vek \lambda$.
        \item 
        Repeat until $\varphi \le 0$:
        \begin{enumerate}[(a), nosep]
        \item
        $s_{\mathrm{old}} \coloneqq s$,
        \item
        $s_{\mathrm{leave}} \coloneqq \min_j \{\mu_j: \mu_j \in (0,1]\}$,
        \item
  $s_{\mathrm{enter}} \coloneqq \min_j \{\mu_j-1: \mu_j >1\}$,
        \item 
        $s \coloneqq s + \min\{s_{\mathrm{leave}},s_{\mathrm{enter}}\}$,
        \item
        $\Vek \mu = \Vek \lambda - s \Vek 1_n$,
        \item  
        $\varphi = \langle \proj_{\mathcal Q} (\Vek \mu), \Vek 1_n \rangle - d$,
        \end{enumerate}
        \item Compute
        \begin{enumerate}[(a), nosep]
        \item
        $k \coloneqq |\{j \in \{1,\ldots,n\}: (\Vek \lambda - s_{\mathrm{old}} \Vek 1_n)_j \in (0,1] \}|$,
        \item 
        $
        \hat t = s_{\mathrm{old}} + \tfrac{1}{k} \varphi(s_{\mathrm{old}}).
        $
        \end{enumerate}
        \end{enumerate}
          {\scshape Output:} $\hat{\Vek \lambda} \coloneqq \proj_{\mathcal{Q}\cap{\mathcal H}}(\Vek \lambda)$.
\end{Algorithm}

Based on the derived proximal mappings,
the primal-dual Algorithm \ref{alg:proj-trunc-cube}
to solve r\pers{Reaper} \eqref{eq:reaper-matrix-form}
can be
specified in matrix form as follows.

\begin{Algorithm}[Primal-dual \mbox{r}\pers{reaper}]   \label{alg:primal-dual-reaper}
  {\scshape Input}: $\Mat X \in \BR^{n, N}$, $d \in \BN$, $\alpha > 0$, and
  $\sigma ,\tau > 0$ with $\sigma \tau < 1/\|\Mat X\|_2^2$, and $\theta \in [0,1)$.\\
  {\scshape Initiation}:
    $\Mat P^{(0)} = \bar{\Mat P}^{(0)}\coloneqq \Mat 0_{n, n}$,
    $\Mat Y^{(0)} \coloneqq \Mat 0_{n, N}$.\\
  {\scshape Iteration}:
    \begin{enumerate}[(i), nosep]
    \item {\scshape Dual update:}
      $\Mat Y^{(r+1)} \coloneqq \proj_{\Set B_{2, \infty}} \bigl( \Vek
      Y^{(r)} + \sigma \bigl(\Op X \bigl(\bar{\Mat P}^{(r)} \bigr) -
      \Mat X \bigr) \bigr)$.
    \item {\scshape Primal update:} 
      \begin{enumerate}[(a)]
      \item 
        $\Mat U \diag(\Vek \lambda) \, \Mat U^\tT \coloneqq \Vek P^{(r)}
        - \tau \Op X^* \bigl(\Vek Y^{(r+1)} \bigr)$,
      \item $\hat{\Vek \lambda} \coloneqq \proj_{_{\mathcal{Q}\cap{\mathcal H}}}(\Vek \lambda -
        \tau \alpha \Vek 1_n)$ \quad (\thref{alg:proj-trunc-cube}),
      \item $\Mat P^{(r+1)} \coloneqq \Mat U \diag(\hat{\Vek \lambda}) \,
        \Mat U^*$.
      \end{enumerate}
    \item {\scshape Extrapolation:}
      $ \bar{\Mat P}^{(r+1)} 
                        \coloneqq (1 + \theta) \, \Mat
      P^{(r+1)} - \theta \, \Mat P^{(r)}$.
    \end{enumerate}
   {\scshape Output}: $\hat{\Mat P}$ (Solution of r\pers{reaper} \eqref{eq:reaper-matrix-form}).
 \end{Algorithm}

\section{Matrix-free realization} \label{sec:matrix-free-real}

Solving r\pers{reaper} with the primal-dual 
\thref{alg:primal-dual-reaper} is possible if the
dimension of the surrounding space $\BR^n$ is moderate which is often not the case in image processing tasks. While
the dual variable $\Mat Y \in \BR^{n,N}$ matches the
dimension of the data, the
primal variable $\Mat P$ is in $S(n)$ instead of $\mathbb  R^{n,d}$, $d \ll n$.
How can the primal-dual iteration be realized in the case $n \gg d$
though the primal variable cannot be hold in memory and the required
eigenvalue decomposition cannot be computed in a reasonable amount of
time?

Here the nuclear norm in r\pers{reaper} that promotes low-rank matrices comes to our
aid.  Our main idea to derive a practical
implementation of the primal-dual iteration is thus
based on the assumption that the iterates of the primal variable $\Mat P^{(r)}$
possess the form
\begin{equation}   \label{eq:low-rank-repr}
  \Mat P^{(r)} 
  \coloneqq 
  \sum_{k=1}^{d_r} \lambda_k^{(r)} \, \Vek u_k^{(r)} \, \bigl(\Vek u_k^{(r)}\bigr)^\tT
 \end{equation}
with small rank $d_r$.  
In our simulations, we observed that
the rank is usually around the dimension $d$ of the wanted
low-dimensional subspace.

In order to integrate the matrix-free representation
\eqref{eq:low-rank-repr} into the primal-dual iteration efficiently,
we further require a fast method to compute the eigenvalue
thresholding.  For this, we compute a partial eigenvalue decomposition
using the well-known Lanczos process \cite{Lan50}.  Deriving
matrix-free versions of the forward operator $\Op X$ and its adjoint
$\Op X^*$, we finally introduce a complete matrix-free primal-dual
implementation with respect to $\Mat P^{(r)}$.

\subsection{The thick-restarted Lanczos process} \label{sec:lancz-proc}

One of the most commonly used methods to extract a small set of
eigenvalues and their corresponding eigenvectors of a large symmetric
matrix is the Lanczos method \cite{Lan50}.  The method
builds a partial orthogonal basis first and then uses a
Rayleigh--Ritz projection to extract the wanted eigenpairs
approximately.  If the set of employed basis vectors is increased, the
extracted eigenpairs converge to the eigenpairs of the given matrix
\cite{GV13}.
Since the symmetric matrix whose partial eigenvalue decomposition is
required in the primal-dual method usually is high-dimensional, we
would like to chose the number $k_{\id{max}}$ of basis vectors
within the Lanczos method as small as possible.  To calculate the
dominant $\ell_{\id{fix}}$ eigenpairs with high accuracy
nevertheless, the Lanczos method can be restarted with the dominant
$\ell_{\id{fix}}$ Ritz pairs.  For our purpose, we use the
thick-restart scheme of Wu and Simon \cite{WS00} in Algorithm \ref{alg:lanczos}, whose
details are discussed below.

\begin{Algorithm}[Thick-restarted Lanczos process {\cite[Alg~3]{WS00}}]
  \label{alg:lanczos}
  {\scshape Input}: $\Mat P \in S(n)$, $k_{\id{max}} >
  \ell_{\id{fix}} > 0$, $\delta > 0$.
  \begin{enumerate}[(i), nosep]
  \item Choose a unit vector $\Vek r_0 \in \BR^n$.  Set
    $\ell \coloneqq 0$.
  \item Lanczos process:
  
    \begin{minipage}[t]{0.5\linewidth}
      {\scshape 1. Initiation}:
      \begin{enumerate}[(a), nosep]
      \item $\Vek e_{\ell + 1} \coloneqq \Vek r_\ell / \pNormn{\Vek
          r_\ell}_2$,
      \item $\Vek q \coloneqq \Mat P \Vek e_{\ell + 1}$,
      \item
        $\beta_{\ell + 1} \coloneqq \iProdn{\Vek q}{\Vek e_{\ell +
            1}}$,
      \item $\Vek r_{\ell + 1} \coloneqq \Vek q - \beta_{\ell + 1}
        \Vek e_{\ell + 1} - \sum_{k=1}^\ell \rho_k \Vek e_k$,
      \item $\gamma_{\ell + 1} \coloneqq \pNormn{\Vek r_{\ell + 1}}$.
      \end{enumerate}
    \end{minipage}
    \begin{minipage}[t]{0.45\linewidth}
      {\scshape 2. Interation} ($k = \ell + 2, \dots, k_{\id{max}}$):
      \begin{enumerate}[(a), nosep]
      \item $\Vek e_k \coloneqq \Vek r_{k-1} / \gamma_{k-1}$,
      \item $\Vek q \coloneqq \Mat P \Vek e_k$,
      \item
        $\beta_k \coloneqq \iProdn{\Vek q}{\Vek e_k}$,
      \item $\Vek r_k \coloneqq \Vek q - \beta_k
        \Vek e_k - \gamma_{k-1} \Vek e_{k-1}$,
      \item $\gamma_k \coloneqq \pNormn{\Vek r_k}$.
      \end{enumerate}
    \end{minipage}
    \vspace{3pt}
  \item Compute the eigenvalue decomposition
    $\Mat T = \Mat Y \Mat \Lambda \Mat Y^\tT$ of $\Mat T$ in \eqref{eq:tridiag-T}.  Set
    $\Mat U \coloneqq \Mat E \Mat Y$.
  \item If $\gamma_{k_{\id{max}}} \absn{y_{k_{\id{max}},k}}
    \le \delta \pNormn{\Mat P}$ for $k = 1, \dots, \ell_{\id{fix}}$,\\
    then return
    $\Mat U \coloneqq [\Vek u_1 | \dots | \Vek u_{\ell_{\id{fix}}}]$
    and $\Mat \Lambda \coloneqq \diag(\lambda_1, \dots,
    \lambda_{\ell_{\id{fix}}})$.\\
    Otherwise, set $\ell \coloneqq \ell_{\id{fix}}$, $\Vek r_\ell
    \coloneqq \Vek r_{k_{\id{max}}}$, and continue
    with (ii).
  \end{enumerate}
  {\scshape Output}: $\Mat U \in \BR^{n \times \ell_{\id{fix}}}$,
  $\Mat \Lambda \in \BR^{\ell_{\id{fix}} \times
    \ell_{\id{fix}}}$
  with $\Mat U^\tT \Mat P \Mat U = \Mat \Lambda$.
\end{Algorithm}

\begin{Remark}
  \label{rem:lanczos:1}
  Although the Lanczos process computes an orthogonal basis $\Vek
  e_1, \dots, \Vek e_{k_{\id{max}}}$, the orthogonality is usually
  lost because of the floating-point arithmetic.  In order to
  re-establish the orthogonality, we therefore have to orthogonalize
  the newly computed $\Vek e_k$ with the previous basis vectors, which
  can be achieved by the Gram--Schmidt procedure.  More
  sophisticated re-orthogonalization strategies are discussed in
  \cite{WS00}.  \qed
\end{Remark}

\begin{Remark}
  \label{rem:lanczos:2}
  During the Lanczos process, the norm of the residual $\gamma_k$
  could become zero.  In this case, we can stop the process, reduce
  $k_{\id{max}}$ to the current $k$, and proceed with step (iii)
  and (iv).  Then the computed basis $\Vek e_1, \dots, \Vek e_k$  spans
  an invariant subspace of $\Mat P$ such that the eigenpairs in
  $\Mat U$ and $\Mat \Lambda$ become exact, see \cite{GV13}.  \qed
\end{Remark}

The heart of the Lanczos method in \thref{alg:lanczos} is the
construction of an orthonormal matrix
$\Mat E \coloneqq [\Vek e_1 | \dots | \Vek e_{k_{\id{max}}}] \in
\BR^{n \times k_{\id{max}}}$
such that $\Mat T \coloneqq \Mat E^\tT \Mat P \Mat E$ becomes
tridiagonal, see \eqref{eq:tridiag-T} with $\ell = 0$ below.  Using
the eigenvalue decomposition $\Mat T = \Mat Y \Mat \Lambda \Mat Y^\tT$,
we then compute the Ritz pairs $(\lambda_k, \Vek u_k)$, where
$\Vek u_k$ are the columns of
$\Mat U \coloneqq [\Vek u_1 | \dots | \Vek u_{k_{\id{max}}}]$ and
$\lambda_k$ the eigenvalues in $\Mat \Lambda$.  In the next iteration,
we chose $\ell_{\id{fix}}$ Ritz pairs corresponding to the
absolute leading Ritz values denoted by
$(\breve \lambda_1, \breve{\Vek u}_1), \dots, (\breve
\lambda_{\ell_{\id{fix}}}, \breve{\Vek u}_{\ell_{\id{fix}}})$
and restart the Lanczos process.  Thereby, the chosen Ritz
vectors are extended to an orthogonal basis
$\Mat E \coloneqq [\breve{\Vek u}_1 | \dots | \breve{\Vek
  u}_{\ell_{\id{fix}}} | \Vek e_{\ell_{\id{fix}}+1} | \dots |
\Vek e_{k_{\id{max}}}]$ fulfilling
\begin{equation}
  \label{eq:tridiag-T}
  \Mat E^* \Mat P \Mat E
  = \Mat T =
  \begin{bmatrix}
    \breve \lambda_1 & & & \rho_1 \\
    & \ddots & & \vdots \\
    & & \breve \lambda_\ell & \rho_\ell \\
    \rho_1 & \cdots & \rho_\ell & \beta_{\ell+1} & \gamma_{\ell+1} \\
    & & & \gamma_{\ell+1} & \beta_{\ell+1}  & \ddots \\
    & & & & \ddots  & \ddots & \gamma_{k_{\id{max}}-1} \\
    & & & & & \gamma_{k_{\id{max}}-1} & \beta_{k_{\id{max}}} \\
  \end{bmatrix},
\end{equation}
where
$\rho_k \coloneqq \breve \gamma_{k_{\id{max}}} \breve
y_{k_{\id{max}}, k}$
with $\breve \gamma_{k_{\id{max}}}$ and
$\breve y_{k_{\id{max}}, k}$ originating from the last iteration,
see \cite{WS00}.

The stopping criteria of the thick-restarted Lanczos process is
here deduced from the fact that the chosen Ritz pairs fulfil the
equation
\begin{equation*}
  \Mat P \, \breve{\Vek u}_k 
  = 
  \breve \lambda_k \, \breve{\Vek u}_k 
  + \breve y_{k_{\id{max}}, k} \, \breve{\Vek r}_{k_{\id{max}}},
\end{equation*}
where $\breve{\Vek r}_{k_{\id{max}}}$ is the last residuum vector
of the previous iteration \cite{WS00}.  Consequently, the absolute
error of the chosen Ritz pairs is given by
\begin{equation*}
  \pNormn{\Mat P \, \breve{\Vek u}_k - 
  \breve \lambda_k \, \breve{\Vek u}_k}_2 
  = \absn{\breve y_{k_{\id{max}}, k}} \,  \pNormn{\breve{\Vek
    r}_{k_{\id{max}}}}_2
  =\breve \gamma_{k_{\id{max}}} \, \absn{\breve y_{k_{\id{max}}, k}}.
\end{equation*}
Usually, the absolute value of the leading Ritz value is a good
approximation of the required spectral norm $\pNormn{\Mat P}$ to
estimate the current relative error.

\subsection{Matrix-free primal update}\label{sec:matr-free-primal}

The thick-restarted Lanczos method allow us to compute the leading
absolute eigenvalues and their corresponding eigenvectors in a
matrix-free manner using only the action of the considered matrix.  In
our primal-dual method for r\pers{reaper}, we need the
action of $\Mat P^{(r)} - \tau  \Op X^*(\Mat Y^{(r+1)})$.  
Incorporating the low-rank representation \eqref{eq:low-rank-repr}, we
see that this can be rewritten as
\begin{equation*}
  \Vek e \in \BR^n
  \mapsto
  \biggl\{
  \sum_{k=1}^{d_r} \lambda_k^{(r)} \, \iProdb{\Vek e}{\Vek u_k^{(r)}}
  \, \Vek u_k^{(r)} \biggr\}
  - \frac{\tau}{2} \, \biggl\{ \Mat Y^{(r+1)} \, \bigl[\Mat X^\tT \Vek e
  \bigr]
  + \Mat X \, \bigl[\bigl( \Mat Y^{(r+1)} \bigr)^\tT \Vek e\bigr] \biggr\}.
\end{equation*}
For the evaluation of the primal proximal mapping, we first compute the
eigenvalue decomposition of
$\Mat P^{(r)} - \tau \Op X^*(\Mat Y^{(r+1)})$, 
next shift the
eigenvalues, and finally project them to the truncated hypercube
$\mathcal{Q} \cap \mathcal{H}$, see \thref{alg:primal-dual-reaper}.  Since the projection
onto $\mathcal{Q} \cap \mathcal{H}$ is independent of negative eigenvalues, see note
after Proposition \ref{thm:proj_hyperqube}, it is thus sufficient to compute only the
eigenpairs with eigenvalue larger than $\alpha \tau$.

For the numerical implementation, we compute the relevant eigenpairs
with the thick-restarted Lanczos method.  In the course of this, we
are confronted with the issue that we actually do not know how many
eigenpairs has to be computed.  To reduce the overhead of
\thref{alg:lanczos} as much as possible, the parameters
$\ell_{\id{fix}}$ and $k_{\id{max}}$ can be easily adapted
between the restarts.  Further, the computation of strongly negative
eigenvalues can be avoided by an eigenvalue shift, \ie\ actually
compute the eigenpairs of
$\Mat P^{(r)} - \tau  \Op X^*(\Mat Y^{(r+1)}) + \nu \Mat I$ with
$\mu \ge 0$, where the required action has the form
\begin{equation}
  \label{eq:act-prim-update}
  \Vek e \in \BR^n
  \mapsto
  \biggl\{
  \sum_{k=1}^{d_r} \lambda_k^{(r)} \, \iProdb{\Vek e}{\Vek u_k^{(r)}}
  \, \Vek u_k^{(r)} \biggr\}
  - \frac{\tau}{2} \, \biggl\{ \Mat Y^{(r+1)} \, \bigl[\Mat X^\tT \Vek e
  \bigr]
  + \Mat X \, \bigl[\bigl( \Mat Y^{(r+1)} \bigr)^\tT \Vek e\bigr]
  \biggr\}
  + \nu \,\Mat e.
\end{equation}
Essentially, we may thus implement the primal proximation in the
following manner.

\begin{Algorithm}[Matrix-free primal proximation]
  \label{alg:primal-prox}
  {\scshape Input}: $\Mat P^{(r)} \in \mathcal S(n)$,
  $\Mat Y^{(r+1)} \in \BR^{n, N}$, $d > 0$, $\tau > 0$,
  $\alpha > 0$.
  \begin{enumerate}[(i), nosep]
  \item {\scshape Thick-restarted Lanczos method}:\\
    Setting $\nu \coloneqq 0$,
    $\ell_{\id{fix}} \coloneqq \rank(\Mat P^{(r)})$,
    $k_{\id{max}} \coloneqq \min\{2 \ell_{\id{fix}}, n\}$, run
    \thref{alg:lanczos} with action \eqref{eq:act-prim-update}.
    Between restarts, check convergence and update parameters:
    \begin{enumerate}[(a), nosep]
    \item If
      $\gamma_{k_{\id{max}}} \absn{y_{k_{\id{max}},k}} \le
      \delta \, \pNormn{\Mat P^{(r)} - \tau  \Op X^*(\Mat Y^{(r+1)})
        + \nu \Mat I}$
      for $k = 1, \dots, m+1$, and\\
      if $\lambda_1 \ge \dots \ge \lambda_m \ge \alpha \tau + \nu >
      \lambda_{m+1}$,\\
      then return $\Mat U \coloneqq [\Vek u_1 | \dots | \Vek u_m]$ and
      $\Mat \Lambda \coloneqq \diag(\lambda_1 - \nu, \dots, \lambda_m
      - \nu)$.
    \item If $\lambda_{\ell_{\id{fix}}} > \alpha \tau + \nu$, then
      increase $\ell_{\id{fix}}$, $k_{\id{max}}$ so that
      $\ell_{\id{fix}} < k_{\id{max}} \le n$.
    \item Set
      $\xi \coloneqq \max\{[\lambda_1]_-, \dots,
      [\lambda_{k_{\id{max}}}]_-\}$ and $\nu \coloneqq \nu + \xi$.
      \\
      Restart with $(\lambda_k + \xi, \Vek u_k)$, $k=1, \dots, \ell_{\id{fix}}$.
    \end{enumerate}
  \item  {\scshape Projection onto} $\mathcal{Q} \cap \mathcal{H}$:\\
    Run \thref{alg:proj-trunc-cube} on $\Vek \lambda \coloneqq
    (\lambda_1 - \alpha \tau, \dots, \lambda_m - \alpha \tau, 0,
    \dots, 0)^\tT \in \BR^n$\\
    to get $\hat {\Vek \lambda} \coloneqq
    \proj_{\mathcal{Q} \cap \mathcal{H}}(\Vek \lambda)$.
  \item  {\scshape New low-rank representation}:\\
    Determine $d_{r+1} \coloneqq \max \{k : \hat \lambda_k > 0\}$ and
    return $\Mat P^{(r+1)} \coloneqq \sum_{k=1}^{d_{r+1}} \hat \lambda_k \, \Vek
    u_k \Vek u_k^\tT$.
  \end{enumerate}
  {\scshape Output}:
  $\Mat P^{(r+1)} \coloneqq \sum_{k=1}^{d_{r+1}} \lambda^{(r+1)}_k
  \, \Vek u_k^{(r+1)} \bigl( \Vek  u_k^{(r+1)} \bigr)^\tT$.
\end{Algorithm}

\begin{Remark}
  If the matrix $\Mat P^{(r)} - \tau \Op X^*(\Mat Y^{(r+1)})$ does
  not possess any eigenvalues greater than $\alpha \tau$, then the
  Lanczos process stops in step (i.a) with $m=0$.  Since the projection
  to the truncated hypercube is then the zero vector again, the new
  iteration $\Mat P^{(r+1)}$ can be represented by an empty low-rank
  representation, \ie\ $d_{r+1}=0$.  \qed
\end{Remark}

\subsection{Matrix-free dual update} \label{sec:matr-free-dual}

Compared with the primal update, the derivation of the matrix-free
dual update is more straightforward.  First, the matrix 
\begin{equation*}
  \Mat Z
  \coloneqq \Mat Y^{(r)} + \sigma \bigl[ \Op X \bigl((1+\theta) \, \Mat P^{(r)} -
  \theta \, \Mat P^{(r-1)}\bigr) - \Mat X \bigr]
\end{equation*}
is computed, where the over-relaxation
$\bar{\Mat P}^{(r)} \coloneqq (1+\theta) \, \Mat P^{(r)} - \theta \,
\Mat P^{(r-1)}$
is already plugged in.  
The low-rank representations of $\Mat P^{(r)}$
and $\Mat P^{(r-1)}$ similar to \eqref{eq:low-rank-repr} can
efficiently incorporated by calculating the matrix
$\Mat Z \coloneqq [\Vek z_1 | \dots | \Vek z_N]$ column by column.
This way of handling the forward operator $\Op X$ nicely matches with
the projection of the columns $\Vek z_k$ to the Euclidean unit
ball in the second step.  Writing the matrix
$\Mat Y^{(r)} \coloneqq [\Vek y_1^{(r)}| \dots | \Vek y_N^{(r)}]$
column by column too, we obtain the following numerical method.

\begin{Algorithm}[Matrix-free dual proximation]
  \label{alg:dual-prox}
  {\scshape Input}: $\Mat Y^{(r)} \in \BR^{n,N}$, $\Mat P^{(r)}
  \in S(n)$, $\Mat P^{(r-1)} \in S(n)$, $\sigma > 0$, $\theta \in
  (0,1]$.
  \begin{enumerate}[(i), nosep]
  \item For $k = 1, \dots, N$, compute
    \vspace*{-12pt}
    \begin{align*}
      \Vek z_k 
      \coloneqq
      \Vek y_k^{(r)}
      + \sigma \, (1 + \theta) \,
      &\biggl\{
        \sum_{\ell=1}^{d_{r}} \lambda_\ell^{(r)} \iProdb{\Vek
        x_k}{\Vek u_\ell^{(r)}} \Vek u_\ell^{(r)}
        \biggr\}
      \\[\fsmallskip]
      - \sigma \theta \,
      &\biggl\{
      \sum_{\ell=1}^{d_{r-1}} \lambda_\ell^{(r-1)} \iProdb{\Vek
        x_k}{\Vek u_\ell^{(r-1)}} \Vek u_\ell^{(r-1)}
      \biggr\} - \sigma \Vek x_k.
    \end{align*}
  \item For $k = 1, \dots, N$, compute
    $\Vek z_k \coloneqq \Vek z_k / (1 + [\pNormn{\Vek z_k}_2 - 1]_+)$.
  \item Return $\Mat Y^{(r+1)} \coloneqq [\Vek z_1 | \dots | \Vek z_N]$.
  \end{enumerate}
  {\scshape Output}: $\Mat Y^{(r+1)}$
\end{Algorithm}

\subsection{Matrix-free projection onto the orthoprojectors}
\label{sec:matr-free-ortho-proj}
With the matrix-free implementations of the primal and dual
proximal mappings, we are already able to solve  r\pers{reaper} \eqref{PCA_3} numerically.  
Before summarizing the compound algorithm, we briefly discuss the last needed component
to tackle the robust PCA problem \eqref{PCA_3_1}.
The final step
is to project the solution $\hat{\Mat P}$ of r\pers{reaper}
onto the set of orthoprojectors with rank not larger than $d$:
\begin{equation*}
\mathcal O_d \coloneqq \{ \Mat \Pi \in \mathcal S(n): \Vek \lambda_{\Vek \Pi} \in \mathcal E_d\},
\end{equation*}
where
\begin{equation*}
\mathcal E_d \coloneqq \{ \lambda \in \BR^n: \Vek \lambda \in \{0,1\}^n, \, \langle \Vek \lambda, \Vek 1_n \rangle \le d\}.
\end{equation*}
We may calculate the projection explicitly in the following manner.

\begin{Proposition}[Projection onto the orthoprojectors]
  \label{prop:proj-orthopr}
  For $\Mat P \in S(n)$ with eigenvalue decomposition
  $\Mat P = \Mat U \diag(\Vek \lambda_{\Mat P}) \, \Mat U^\tT$, and
  for every $1 \le p \le \infty$, the projection onto $\mathcal O_d$
  with respect to the Schatten $p$-norm is given by
  \begin{equation*}
    \proj_{\mathcal O_d}(\Mat P) = \Mat U \diag(\proj_{\mathcal E_d}(\Vek \lambda)) \,     \Mat U^\tT.
  \end{equation*}
\end{Proposition}

\begin{Proof}
  The key ingredient to prove this statement is the theorem of
  Lidskii--Mirsky--Wielandt, see for instance \cite{LM99}.  Using this
  theorem to estimate the Schatten $p$-Norm, we obtain
  \begin{equation}    \label{eq:eigen-esti}
    \min_{\Pi \in \mathcal O_d} \|\Mat P - \Mat \Pi \|_{S_p} 
                \ge 
                \min_{\Vek \lambda_{\Mat \Pi }\in \mathcal E_d} \| \Vek \lambda_{\Mat P }   - \Vek \lambda_{\Mat \Pi }  \|_p,
  \end{equation}
  where we have equality if $\Mat \Pi$ has the same eigenvectors as
  $\Mat P$.  Recall that the eigenvalues in $\Vek \lambda_{\Mat P }$
  appear in descending order.  The right-hand side of
  \eqref{eq:eigen-esti} thus becomes minimal if we choose the
  eigenvalues of $\Mat \Pi$ for $k=1,\dots, d$ as
  \begin{align*}
    \hat \lambda_{\Mat \Pi,k}
    & \coloneqq
      \left\{
      \begin{aligned}
        1 \quad& \text{if } \lambda_{\Mat P,k} \ge \tfrac12, \\
        0 \quad& \text{if } \lambda_{\Mat P,k} < \tfrac12, \\
      \end{aligned}
    \right.
   \end{align*}
        and set
        $
    \hat \lambda_{\Mat \Pi,k}
    \coloneqq 0 
        $
        for $k=d+1,\ldots,n$.
  This is exactly the projection onto  $\mathcal E_d$.    \qed
\end{Proof}

Because of the low-rank representation
$\Mat P^{(r)} = \sum_{k=1}^{d_r} \lambda_k \, \Vek u_k^{(r)} (\Vek
u_k^{(r)})^\tT$ of the primal variable, the construction of the
orthoprojector $\hat{\Mat \Pi} \in \mathcal O_d$ is here especially
simple.

\begin{Algorithm}[Matrix-free projection onto orthoprojectors]
  \label{alg:matr-free-ortho}
  {\scshape Input}:
  $\Mat P = \sum_{k=1}^\kappa \lambda_k \, \Vek u_k \Vek u_k^\tT \in \mathcal S(n)$,
  $d \in \BN$.
  \begin{enumerate}[(i), nosep]
  \item {\scshape Projection onto} $\mathcal O_d$:\\
    Determine $s \coloneqq \max \{ k : \lambda_k \ge \nicefrac12, k
    \le \min\{d, \kappa\}\}$.
  \item {\scshape Matrix-free presentation}:\\
    Return $\Mat \Pi = \sum_{k=1}^s \Vek u_k \Vek u_k^\tT$.
    \end{enumerate}
  {\scshape Output}: $\Mat \Pi = \proj_{\mathcal O_d}(\Mat P)$.
\end{Algorithm}

\subsection{Matrix-free robust PCA by r\pers{reaper}}
\label{sec:matr-free-pca}

Combining the matrix-free implementations of the primal and dual
proximal mappings, we finally obtain a primal-dual method to solve 
r\pers{reaper} \eqref{PCA_3} without
evaluating the primal variable $\Mat P^{(r)}$ representing the relaxed
orthoprojector explicitly.

\begin{Algorithm}[Matrix-free robust PCA]
  \label{alg:matrix-free-reaper}
  {\scshape Input}: $\Mat X \in \BR^{n, N}$, $d \in \BN$, $\alpha > 0$, and
  $\sigma ,\tau > 0$ with $\sigma \tau < 1/\|\Mat X\|_2^2$, and $\theta \in (0,1]$.\\
  {\scshape Initiation}:
    $\Mat P^{(0)} = \bar{\Mat P}^{(0)}\coloneqq \Mat 0 \in \BR^{n, n}$,
    $\Mat Y^{(0)} \coloneqq \Mat 0 \in \BR^{n, N}$.\\
  {\scshape Iteration}: 
        \begin{enumerate}[(i), nosep]
     \item {\scshape Dual update:} Compute $\Mat Y^{(r+1)}$ with
      \thref{alg:dual-prox}. 
    \item {\scshape Primal update:} Compute $\Mat P^{(r+1)}$ with
      \thref{alg:primal-prox}.
    \end{enumerate}
 {\scshape Projection}: Compute $\hat{\Mat \Pi}$ with
    \thref{alg:matr-free-ortho}.
  \\
  {\scshape Output}: $\hat{\Mat \Pi}$ (Orthoprojector onto recovered subspace).
\end{Algorithm}

\section{Performance analysis} \label{sec:error-analysis}
Inspired by ideas of Lerman et al. \cite{LCTZ15},
we examine the performance analysis of r\pers{reaper}.
To this end, we assume that the \bq{ideal} subspace $L$ of the given data 
$\Vek x_k \in \mathbb R^n$, $k=1,\ldots,N$ has dimension $d_L \le d$.
As in \cite{LCTZ15}, we determine the best fit of the data 
by two measures: the first one is the distance of the data from the subspace
\begin{equation*}
\Op R_L = \Op R_L (\Mat X)   \coloneqq \| (\Mat I_n - \Mat \Pi_L) \Mat X \|_{2,1},
\end{equation*}
where $\Pi_L$ denotes the orthogonal projector onto $L$.
For the second measure, we assume that the projected data $\{\Mat \Pi \Vek x_k: k=1,\ldots,N \}$, $N \ge d_L$, 
form a \emph{frame} in $L$ meaning that there exist constants $0 < c_L \le C_L < \infty$ such that
\begin{equation*}
 c_L \le   \sum_{k=1}^N   \absn{\iProdn{\Vek u}{\Mat \Pi_L \Vek x_k}}^2  = \sum_{k=1}^N \absn{\iProdn{\Vek u}{\Vek x_k}}^2 \le C_L 
\end{equation*}
for all $\Vek u \in L$ with $|u\|_2 = 1$.
In order to recover the entire subspace $L$, the data have obviously
to cover each direction in $L$ with sufficiently many data points.  This
 well-localization of the data is measured
by the \emph{permeance statistic}
\begin{equation} \label{permeance}
 \Op P_L = \Op P_L(\Mat X)
  \coloneqq
  \min_{\substack{\Vek u \in L\\ \pNormn{\Vek u}=1}} \sum_{k=1}^N
  \absn{\iProdn{\Vek u}{\Vek x_k}}
\end{equation}
which can be seen as $\ell_1$ counterpart of the lower frame bound.
Clearly, $\Op P_L$ becomes large if
all direction in $L$ are uniformly covered by the data.
The lower frame bound and the permeance statistic come into the play in the following lemma,
compare with \cite[Section A2.3]{LCTZ15}.

\begin{Lemma}\label{lower}
Let $\Mat \Pi_L$ be the orthogonal projector onto a subspace $L$ of $\mathbb R^n$ of dimension $d_L$
and $\Vek x_k \in \mathbb R^n$, $k=1,\ldots,N$, $N \ge d_L$ which form the columns of the matrix $\Mat X$.
Then, for any $\Mat A \in \mathbb R^{n,n}$, the following relations hold true:
\begin{align} \label{lower_1}
 \| \Mat A \Mat \Pi_L \Mat X\|_{2,2} 
&\ge 
c_L \| \Mat A \Mat \Pi_L \|_F,
\end{align}
and
\begin{align}\label{lower_2}
\| \Mat A \Mat \Pi_L \Mat X\|_{2,1} 
&\ge 
\Op P_L \| \Mat A \Mat \Pi_L \|_F 
\ge 
\tfrac{1}{ \sqrt{d_L} } \Op P_L \| \Mat A \Mat \Pi_L \|_{\tr}.
\end{align}
\end{Lemma}

\begin{Proof} We restrict our attention to \eqref{lower_2}.
The relation \eqref{lower_1} follows similar lines.
Let $\Mat A \,  \Mat \Pi_L$ have the singular value decomposition
$\Mat A \, \Mat \Pi_L = \Mat U \Mat \Sigma \Mat V^\tT$, where the singular values $\sigma_k$, $k=1,\ldots,n$ 
are in descending order and $\sigma_{d_L+1}= \ldots = \sigma_n = 0$ and
we can arrange $\Mat V$ such that the transpose of the first $d_L$ rows of $\Mat V$ belong to  $L$.
Then it holds
$$\|\Mat A \Mat \Pi_L\|_F^2 = \sum_{k=1}^{d_L} \sigma_k^2.$$ 
Using orthogonality of $\Mat U$ and concavity of the square root function, we obtain
\begin{align*}
\|\Mat A \Mat \Pi_L X\|_{2,1} 
&= \| \Mat U \Mat \Sigma \Mat V^\tT \Mat X \|_{2,1}
= \|  \Mat \Sigma \Mat V^\tT \Mat X \|_{2,1}\\
&= \sum_{k=1}^N \left( \sum_{j=1}^{d_L} \sigma_j^2 \langle \Vek v_j,\Vek x_k \rangle^2 \right)^\frac12 \\
&= \|\Mat A \Mat \Pi_L\|_F \,  \sum_{k=1}^N \left( \sum_{j=1}^{d_L} \frac{\sigma_j^2}{\|\Mat A \Mat \Pi_L\|_F^2}  
\langle v_j,x_k \rangle^2 \right)^\frac12\\
&\ge 
\|\Mat A \Mat \Pi_L\|_F \,  \sum_{k=1}^N \sum_{j=1}^{d_L} \frac{\sigma_j^2}{\|\Mat A \Mat \Pi_L\|_F^2}
 |\langle \Vek v_k,\Vek x_j \rangle|\\
 &= 
\|\Mat A \Mat \Pi_L\|_F \,  \sum_{j=1}^{d_L} \frac{\sigma_j^2}{\|\Mat A \Mat \Pi_L\|_F^2} 
\sum_{k=1}^N  
 |\langle \Vek v_k, \Vek x_j \rangle|\\
&\ge 
\Op P_L \, \|\Mat A \Mat \Pi_L\|_F . 
\end{align*}
The second estimate in \eqref{lower_2} follows by \eqref{mittel}. \hfill \qed
\end{Proof}

Now we can estimate the
reconstruction error of r\pers{reaper}.

\begin{Theorem}  \label{thm:stab-analy}
 Let $\Mat \Pi_L$ be the orthogonal projector onto a subspace $L$ of $\mathbb R^n$ of dimension $d_L$
and $\Vek x_k \in \mathbb R^n$, $k=1,\ldots,N$, $N \ge d_L$ such that their projections onto $L$ form a frame of $L$.
Define $\Op P_L$  by \eqref{permeance} and set $\gamma_L \coloneqq \frac{1}{2\sqrt{2 d_L}}   \Op P_L$.
Let $\hat{\Mat P}$ be the
  solution of \eqref{PCA_3} and $\hat{\Mat \Pi}$ the
  projection of $\hat{\Mat P}$ onto $\mathcal O_d$. Then, for $\alpha \le 2\gamma_L$ the reconstruction
  error is bounded by
  \begin{equation*}
    \pNormn{\hat{\Mat \Pi} - \Mat \Pi_L}_{\tr} \le \frac{8 \Op R_L}{\gamma_L- |\gamma_L - \alpha|}.
      \end{equation*}
\end{Theorem}

\begin{Proof}
Since $\hat {\Mat P}$ is a minimizer of \eqref{PCA_3}, we obtain
\begin{align}
0 \le& 
\|(\Mat I_n - \Mat \Pi_L)\Mat X\|_{2,1} - \|(\Mat I_n - \hat{\Mat P}) \Mat X\|_{2,1} 
+ \alpha (\| \Mat \Pi_L\|_{\tr} - \|\hat {\Mat P}\|_{\tr}) \label{funk}\\
=& 
\Op R_L - \|(\Mat I_n - \hat{\Mat P}) \Mat X\|_{2,1} 
 + \alpha (d_L - \|\hat {\Mat P}\|_{\tr}) \notag\\ 
\le& 
\Op R_L + \|(\Mat I_n - \hat{\Mat P}) (\Mat I_n -\Mat \Pi_L) \Mat X\|_{2,1} -  \|(\Mat I_n - \hat{\Mat P}) \Mat \Pi_L \Mat X\|_{2,1}    
+ \alpha (d_L - \|\hat {\Mat P}\|_{\tr})\notag\\
=&
\Op R_L + \|\big(\Mat I_n - (\hat{\Mat P} - \Mat \Pi_L) \big) (\Mat I_n -\Mat \Pi_L) \Mat X\|_{2,1} 
-  \|(\Mat I_n - \hat{\Mat P}) \Mat \Pi_L \Mat X\|_{2,1} + \alpha (d_L - \|\hat {\Mat P}\|_{\tr})\notag\\
\le&
\left( 2+\|\hat{\Mat P} - \Mat \Pi_L\|_2        \right) \Op R_L 
-  \|(\Mat I_n - \hat{\Mat P}) \Mat \Pi_L \Mat X\|_{2,1} + \alpha (d_L - \|\hat {\Mat P}\|_{\tr}). \label{zerleg}
\end{align}
It remains to estimate 
$\|(\Mat I_n - \hat{\Mat P}) \Mat \Pi_L \Mat X\|_{2,1} = \|(\hat{\Mat P} - \Mat \Pi_L) \Mat \Pi_L \Mat X\|_{2,1}$
from below.
To this end, we decompose $\Mat \Delta := \hat{\Mat P} - \Mat \Pi_L$ as
\begin{align*}
    \Mat \Delta
    &= \underbrace{ \Mat \Pi_{L} \Mat \Delta \Mat \Pi_{L} }_{\Mat \Delta_1}
    + \underbrace{ (\Mat I_n - \Mat \Pi_L) \Mat \Delta \Mat \Pi_{L} }_{\Mat \Delta_2}
    + \underbrace{ \Mat \Pi_{L} \Mat \Delta (\Mat I_n - \Mat \Pi_L) }_{\Mat \Delta_2^\tT}
    + \underbrace{ (\Mat I_n - \Mat \Pi_L) \Mat \Delta (\Mat I_n - \Mat \Pi_L) }_{\Mat \Delta_3}.
 \end{align*}
Since $\Mat 0_{n,n} \preceq \Mat \hat P \preceq \Mat I_n$, we obtain be conjugation with $\Mat \Pi_L$, resp. $I_n - \Mat \Pi_L$
that $\Mat \Delta_1 \preceq \Mat 0_{n,n}$ and $\Mat 0_{n,n} \preceq \Mat \Delta_3$,
so that
$\| \Mat \Delta_1 \|_{\tr} = -\tr{\Mat \Delta_1}$  and $\| \Mat \Delta_3 \|_{\tr} = \tr{\Mat \Delta_3}$.
Then we conclude
\begin{align*}
\tr (\Mat \Delta) &= \tr(\hat{\Mat P}) - d_L  = \tr{\Mat \Delta_1} + 2 \tr{\Mat \Delta_2} + \tr{\Mat \Delta_3}\\
&= \tr{\Mat \Delta_1} + \tr{\Mat \Delta_3}
= \| \Mat \Delta_3 \|_{\tr} - \| \Mat \Delta_1 \|_{\tr},
\end{align*}
which implies
\begin{align*}
\pNormn{\Mat \Delta}_{\tr}
    &\le \pNormn{\Mat \Delta_1}_{\tr} +\pNormn{\Mat \Delta_2}_{\tr}
      + \pNormn{\Mat \Delta_2^\T}_{\tr} + \pNormn{\Mat \Delta_3}_{\tr}
    \\[\fskip]
    &= 2 \pNormn{\Mat\Delta_1}_{\tr} + 2 \pNormn{\Mat\Delta_2}_{\tr} + \tr(\hat{\Mat P})-d_L.  
\end{align*}
Now we can estimate the last summand by \eqref{mittel} and Lemma \ref{lower} as
\begin{align}
\|\Mat \Delta \Mat \Pi_L \Mat X\|_{2,1}
&= 
\|\Mat \Pi_{L} \Mat \Delta \Mat \Pi_L \Mat X + (\Mat I_n - \Mat \Pi_L) \Mat \Delta \Mat \Pi_L \Mat X\|_{2,1}
\notag\\
 &= 
\sum_{k=1}^N \left( \|\Mat\Delta_1 \Vek x_k\|_2^2 + \|\Mat\Delta_2 \Vek x_k\|_2^2\right)^\frac12\notag\\
&\ge
\frac{1}{\sqrt{2}} \sum_{k=1}^N \left( \|\Mat\Delta_1 \Vek x_k\|_2 
+\|\Mat\Delta_2 \Vek x_k\|_2 \right)\notag\\
&\ge
\frac{1}{ \sqrt{2 d_L}} \Op P_L 
\left( \| \Mat\Delta_1\|_{\tr} 
+ 
\|\Mat\Delta_2 \|_{\tr} \right)\notag\\
&\ge 
\frac{1}{2\sqrt{2 d_L}}   \Op P_L   \left( \|\Mat \Delta \|_{\tr} + d_L -  \|\hat{\Mat P}\|_{\tr} \right). 
\label{lower3}
\end{align}
By \eqref{zerleg} and \eqref{lower3}, we obtain
\begin{align}\label{ende}
0 
&\le (2 + \|\hat{\Mat P} - \Mat \Pi_L\|_2 ) \, \Op R_L + (\alpha - \gamma_L) \left(d_L -  \|\hat{\Mat P}\|_{\tr}\right)
-  \gamma_L  \|\hat{\Mat P} - \Mat \Pi_L\|_{\tr}  \\
&\le
(2 + \|\hat{\Mat P} - \Mat \Pi_L\|_2 ) \, \Op R_L + |\alpha - \gamma_L| \left|d_L -  \|\hat{\Mat P}\|_{\tr}\right|
-  \gamma_L  \|\hat{\Mat P} - \Mat \Pi_L\|_{\tr} .\notag
\end{align}
Using that by the triangular inequality
$$
\|\hat{\Mat P} - \Mat \Pi_L\|_{\tr} \ge \left| d_L - \| \hat{\Mat P}\|_{\tr} \right|
$$
we get
\begin{align*}
0 &\le 
(2 + \|\hat{\Mat P} - \Mat \Pi_L\|_2 ) \, \Op R_L + \left(|\alpha - \gamma_L| - \gamma_L \right) \ \|\hat{\Mat P} - \Mat \Pi_L\|_{\tr} .
\end{align*}
Now we can use the estimate $\|\hat{\Mat P} - \Mat \Pi_L\|_2 \le 2$ 
to get 
$$
\|\hat{\Mat P} - \Mat \Pi_L\|_{\tr} \le \frac{4 \Op R_L}{\gamma_L - |\gamma_l - \alpha|}
$$
if $\alpha < 2 \gamma_L$.
The final assertion follows by
\begin{equation*}
\|\hat{\Mat \Pi} - \hat{\Mat \Pi}_L\|_{\tr} \le \|\hat{\Mat \Pi} - \hat{\Mat P}\|_{\tr}
+\|\hat{\Mat P} - \hat{\Mat \Pi}_L\|_{\tr} \le 2 \|\hat{\Mat P} -
\hat{\Mat \Pi}_L\|_{\tr}.
\tag*{\qed}
\end{equation*}
\end{Proof}

\begin{Remark}
Using
$\|\hat{\Mat P} - \Mat \Pi_L\|_2 \le \|\hat{\Mat P} - \Mat \Pi_L\|_{\tr}$
we could alternatively estimate
$$
\|\hat{\Mat P} - \Mat \Pi_L\|_{\tr} \le \frac{2 \Op R_L}{\gamma_L - |\gamma_l - \alpha| - \Op R_L}
$$
if $\gamma_L - |\gamma_l - \alpha| - \Op R_L > 0$.

Further, if $\Vek x_k \in L$, $k=1,\ldots,N$, then, we have by \eqref{funk} that
$$
0 \le - \|(\Mat I_n - \hat{\Mat P}) \Mat X\|_{2,1} 
+ \alpha ( d_L - \|\hat {\Mat P}\|_{\tr}) \label{funk_1}
$$
so that $d_L \ge \|\hat {\Mat P}\|_{\tr}$.
On the other hand, we know by \eqref{ende} that
$$
\|\hat{\Mat P} - \Mat \Pi_L\|_{\tr}  \le \frac{\alpha - \gamma_L}{\gamma_L} \left(d_L - \tr(\hat{\Mat P}) \right)
$$ 
which is not possible if $d_L > \tr(\hat{\Mat P})$ and $\alpha \le \gamma_L$,  so that for such an $\alpha$ we would get 
$\tr(\hat{\Mat P} ) = d_L$.
 \qed
\end{Remark}

\section{Incorporating the offset} \label{sec:incorporating-offset}
So far, we have assumed that the offset $\Vek b$ in the robust PCA problem is given, so that we can just search
for a low dimensional linear subspace which represents the data well.
While in the classical PCA \eqref{PCA_1} the affine subspace always passes trough the mean value
\eqref{mean} of the data, it is not clear which value $\Vek b$ must be chosen in order to minimize
\begin{equation} \label{pca_offset}
 E(\Mat A,\Vek b) \coloneqq      \sum_{i=1}^N \|(\Mat A \Mat A^\tT - \Mat I_n )(\Vek b-\Vek x_i)\|_2 \quad \mathrm{subject \; to} \quad \Mat A^\tT \Mat A = \Mat I_d.
\end{equation}
Clearly, if $(\hat {\Mat A}, \hat{\Vek b})$ is a minimizer of $E$, then, 
for every $\Vek b \in \ran(\hat{\Mat A})$, also $(\hat {\Mat A}, \hat{\Vek b} + \Vek b)$ 
is a minimizer.

A common choice for the offset is the \emph{geometric median} of the data points defined by
\begin{equation*}
\hat{\Vek b} \coloneqq \argmin_{\Vek b \in \mathbb R^n} \sum_{k=1}^N \|\Vek b - \Vek x_k\|_2,
\end{equation*}
which can be computed e.g. by the Weiszfeld algorithm and its
generalizations, see, e.g. \cite{BS2015,OL1978,SST2012,We37}.  Other
choices arising, e.g. from Tylor's M-estimator or other robust
statistical approaches \cite{KTV94,LNNS2019,LM2018,MVD97,Tyl87}, were
proposed in the literature.  However, they do in general not minimize
\eqref{pca_offset} as the following example from \cite{NNSS19b} shows:
given three points in $\mathbb R^n$ which form a triangle with angles
smaller than 120 degrees, the geometric median is the point in the
inner of the triangle from which the points can be seen under an angle
of 120 degrees.  In contrast, the line ($d=1$) having smallest
distance from the three points is the one which passes through those
two points with the largest distance from each other.

In the following, we show that there always exists an optimal
hyperplane of dimension $d=n-1$ in $\mathbb R^n$ determined by a
minimizer of $E$ in \eqref{pca_offset} that
contains at least $n$ data points. Further, if
the number $N$ of data points is
odd, then \emph{every} optimal hyperplane contains at
least $n$ data points.  Recall that the hyperplane spanned by the
columns of
$\Mat A = (\Vek a_1|\ldots|\Vek a_{n-1}) \in \mathrm{St}(n,n-1)$ with
offset $\Vek b$ is given by
$$
\{\Vek x = \Mat A \Vek t + \Vek b: \Vek t \in \BR^{n-1} \} 
=
\{\Vek x \in \BR^n: \langle \Vek a^\perp , \Vek x \rangle = \beta\},
$$
where $\Vek a^\perp \perp \Vek a_i$, $i=1,\ldots,n-1$ is a unit normal vector of the hyperplane, which is uniquely determined  up to its sign
and $\langle \Vek a^\perp,\Vek b \rangle = - \beta$.

The following lemma describes the placement of the data points with respect to the halfspaces 
determined by a minimizing hyperplane.

\begin{Lemma}\label{lem:one_point}
Let  $\Vek x_k \in \BR^n$, $k=1,\ldots,N$. Let $(\hat{\Mat A},\hat{\Vek b})$ be a minimizer of $E$
and $N=M + M_+ + M_-$ with
\begin{align*}
M \coloneqq &|\{\Vek x_k: \langle \hat{\Vek a}^\perp, \Vek x_i \rangle = \hat \beta\}|,\\
M_+ \coloneqq &|\{\Vek x_k: \langle \hat{\Vek a}^\perp, \Vek x_i \rangle > \hat \beta\}|,\\
M_- \coloneqq &|\{\Vek x_k: \langle \hat{\Vek a}^\perp, \Vek x_i \rangle < \hat \beta\}|.
\end{align*}
Then it holds $|M_+ - M_-| \le M$.
In particular, it holds $M \ge 1$ if $N$ is odd.
Also for even $N$ there exists a minimizing hyperplane with $\hat{\Vek b} = \Vek x_k$ for some $k \in \{1,\ldots,N\}$.
\end{Lemma}

\begin{Proof}
W.l.o.g. assume that $M_+ \ge M_-$. If $M_+ = 0$, then all data points are in the mimimizing hyperplane and we are done.
Otherwise, the value $\varepsilon \coloneqq \min\{ \langle \hat{\Vek a}^\perp, \Vek x_k \rangle > 0: k=1,\ldots,N \}$ is positive,
and we consider the shifted hyperplane
$\{\Vek x \in \BR^n: \langle \Vek a^\perp , \Vek x \rangle = \beta + \varepsilon\}$,
which contains at least one data point.
The sum of the distances of the data points from this hyperplane is 
$$
E(\hat{\Mat A},\hat{\Vek b}) - \varepsilon M_+ + \varepsilon (M+M_-).
$$
Since $(\hat{\Mat A},\hat{\Vek b})$ is a minimizer of $E$ this implies that
$M_+ \le M + M_-$.
If $M = 0$, then $M_- = M_+$ and the shifted hyperplane is also minimizing.
However, this case cannot appear for odd $N$ so hat $M \ge 1$ for odd $N$.
This finishes the proof.  \qed
\end{Proof}

\begin{Theorem} \label{thm:offset}
Let  $\Vek x_\ell \in \BR^n$, $\ell=1,\ldots,N$. Then there exists a minimizer $(\hat{\Mat A},\hat{\Vek b})$ of $E$
such that the corresponding minimizing hyperplane contains at least $n$ data points.
If $N$ is odd, every minimizing hyperplane contains at least $n$ data points.
\end{Theorem}

\begin{Proof}
  By Lemma \ref{lem:one_point} there exists a data point $\Vek x_\ell$
  such that $(\hat {\Mat A}, \hat {\Vek b})$ with
  $\hat {\Vek b} = \Vek x_\ell$ is a minimizer of $E$ and for odd $N$
  every minimizing hyperplane passes through a data point.  Let
  $\hat{\Vek a}^\perp$ be a unit vector orthogonal to the columns of
  $\hat {\Mat A}$.  Set $\Vek y_k \coloneqq \Vek x_k - \Vek x_\ell$,
  $k=1,\ldots,N$.  
        Next, we show: if the subspace normal to
    $\hat{\Vek a}^\perp$ contains $M$ linearly independent vectors
    $\Vek y_1, \dots, \Vek y_M$ with $0 \le M \le n-2$, then exactly
    one of the following situations apply.  
                (i) The remaining vektors
    $\Vek y_k$ with $k=M+1, \dots, N$ are linearly dependent from the
    first $M$ vectors and thus in the same linear subspace 
                $\mathrm{span} \{\Vek y_k: k=1,\ldots,M\}$. 
                (ii) We
    find a further independent vector, say $\Vek y_{M+1}$, contained
    in the subspace normal to $\hat{\Vek a}^\perp$ such that we can
    increase $M$ to $M+1$.
  Repeating this argumentation until
  $M = n-2$, we are done since $\Vek y_k + \Vek x_\ell$,
  $k=1,\ldots,n-1$ and $\Vek x_\ell$ itself are in the subspace
  corresponding to $(\hat{\Mat A}, \hat{\Vek b})$.

  Because the vectors $\Vek y_k$ with
    $k=1,\ldots,M$ are independent and are contained in the subspace
    normal to $\hat{\Vek a}^\perp$ by assumption, there exists a matrix
  $$\Mat B \coloneqq (\Vek u_1|\ldots|\Vek u_M| \Vek u_{M+1}|\ldots| \Vek u_{n-1}) \in \mathrm{St}(n,n-1),$$ 
 with $\ran(\Mat A) = \ran(\Mat B)$, whose first
    columns have the same span as $\Vek y_1, \dots, \Vek y_M$,
    \ie\
  $$
  \mathrm{span}\{\Vek u_k:k=1,\ldots,M\} = \mathrm{span}\{\Vek y_k:k=1,\ldots,M\}.
  $$
 This especially implies $\Vek y_k \perp \Vek u_\ell$
    for $k = 1, \dots, M$ and $\ell = M+1, \dots, n-1$.  Note that
  the normal unit vector of $\ran(\Mat B)$ is also
  $\hat{\Vek a}^\perp$, and that the objectives
    coincides, \ie\
  \begin{equation} \label{equiv_1}
    \sum_{k=1}^N \|( \hat{\Mat A}  \hat{\Mat A}^\tT - I_n) \Vek y_k \|_2 = \sum_{k=1}^N \|( \Mat B  \Mat B^\tT - I_n) \Vek y_k \|_2. 
  \end{equation}
Now, let the matrix-valued function
  $\phi_{\Mat B}: [-\pi,\pi) \rightarrow \mathrm{St}(n,n-1)$ be
  defined by
  $$
  \phi_{\Mat B}(\alpha) := \Mat Q_{\Mat B} \, \Mat R(\alpha) \, \Mat C,
  $$
  where the three building factors are given by
  $$
  \Mat Q_{\Mat B} \coloneqq (\Mat B|\hat{\Vek a}^\perp), \quad
  \Mat R(\alpha) \coloneqq
  \left(
    \begin{array}{ccc}
      \Mat I_{n-2}     &\Vek 0_{n-2}&\Vek 0_{n-2}\\[0.5ex]
      \Vek 0_{n-2}^\tT &\cos(\alpha)&\sin(\alpha)\\[0.5ex]
      \Vek 0_{n-2}^\tT &-\sin(\alpha)&\cos(\alpha)
    \end{array}
  \right),
  \quad
  \Mat C \coloneqq
  \begin{pmatrix}                                                                             
    \Mat I_{n-1}\\
    \Vek 0_{n-1}^\tT
  \end{pmatrix}.
  $$
  Figuratively, the function $\phi_{\Mat B}$ takes the
  orthonormal columns of $\Mat B$ and rotates the last vektor $\Vek
  u_{n-1}$ by the angle $\alpha$ in the plane spanned by $\Vek
  u_{n-1}$ and $\hat{\Vek a}^\perp$.  
  Clearly, we have $\phi_{\Mat B}(0) = \Mat B$.
  Due to \eqref{equiv_1}, the function 
  $$
  F(\alpha) = \sum_{k=1}^N \|( \phi_{\Mat B}(\alpha) \phi_{\Mat B}(\alpha)^\tT- \Mat I_n) \Vek y_k \|_2 = \sum_{k=1}^N f_k(\alpha) 
  $$
  has moreover a minimum in $\alpha = 0$. For the
  summands of $F$, we obtain
  \begin{align*}
    f_k(\alpha) 
    &= 
      \|( 
      \phi_{\Mat B}(\alpha) \phi_{\Mat B}(\alpha)^\tT - \Mat I_n 
      ) \, \Vek y_k \|_2\\
    &= 
      \| ( \Mat Q_{\Mat B} \, \Mat R(\alpha) \, \Mat C \, \Mat C^\tT \, \Mat R(\alpha)^\tT \,  \Mat Q_{\Mat B}^\tT  
      - \Mat I_n ) \,  \Vek y_k \|_2\\
    &=  
      \|(
      \Mat C \, \Mat C^\tT \, \Mat R(\alpha)^\tT \,  \Mat Q_{\Mat B}^\tT - \Mat R(\alpha)^\tT \, \Mat Q_{\Mat B}^\tT 
      ) \, \Vek y_k \|_2\\
    &= 
      \|(\Mat C \, \Mat C^\tT - \Mat I_n) \Mat R(\alpha)^\tT \,  \Mat Q_{\Mat B}^\tT \Vek y_k \|_2\\
    &= 
      \left| \sin(\alpha) \langle \Vek u_{n-1},   \Vek y_k \rangle    + \cos(\alpha) \langle \hat{\Vek a}^\perp,   \Vek y_k \rangle \right|,
  \end{align*}
  since $\Mat Q_{\Mat B}$ and $\Mat R(\alpha)$ are
    orthogonal by construction.
  Hence, we get
  $$
  F(\alpha) = \sum_{k=M+1}^N  \left| \sin(\alpha) \langle \Vek u_{n-1},   \Vek y_k \rangle    + \cos(\alpha) \langle \hat{\Vek a}^\perp,   \Vek y_k \rangle \right|.
  $$
 Here the first $M$ summands vanish because of the
    mentioned orthogonality $\Vek y_k \perp \Vek u_{n-1}$ and $\Vek
    y_k \perp \hat{\Vek a}^\perp$ for $k=1, \dots, M$.
        
        If all remaining given points $\Vek y_k$ with
  $k=M+1,\ldots,N$ are in
  $\mathrm{span}\{\Vek y_1, \ldots, \Vek y_M\}$, then
  the corresponding remaining summands of $F(\alpha)$
  become zero too, and the first situation (i) applies; so we are done.

  If this is not the case, consider only those
  $\Vek y_k$ with $k=M+1,\ldots,N$ that are linearly independent of the
  $\Vek y_k$, $k=1,\ldots,M$. Let us denote the corresponding non-empty
  index set by $\mathcal I$.  Assume that there exists a
  $k \in \mathcal I$ such that $f_k$ is not differentiable in
  $\alpha = 0$.  This is only possible if the argument
  of the absolute value vanishes implying
  $$
  f_k(0) = \absn{\langle \hat{\Vek a}^\perp,   \Vek y_k \rangle} = 0.
  $$
  Thus, the vector $\Vek y_k$ is in the
  subspace spanned by the columns of $\hat{\Mat A}$, and we are done.  Otherwise,
  if $f_k(0) \not = 0$ for all $k\in \mathcal I$, then it is
  differentiable in $\alpha = 0$ and, by straightforward
  differentiation, we obtain
  $$
  f_k''(0) = - f_k(0) < 0.
  $$
  But then $\alpha = 0$ cannot be a minimizer of $F$ which is a contradiction.
  Hence, this case cannot occur and the proof is complete.  \qed
\end{Proof}

If the target dimension $d$ of the minimizing
subspace is strictly less than $n-1$,
then it does not have to contain any data point as the following
example shows.

\begin{Counterexample}[Lower-dimensional subspace approximation]
  Initially, we consider the approximation of some given
    points in $\BR^3$ by an one-dimensional subspace -- a line.  More
    precisely, for a fixed $T \gg 1$, we consider the six given points
  \begin{equation*}
    (\cos(\phi),
    \sin(\phi),
    \pm T)^\T
    \qquad\text{with}\qquad
    \phi \in \{0, \nicefrac{2\uppi}{3}, \nicefrac{4\uppi}{3}\}.
  \end{equation*}
  We thus have two well-separated clusters
  around $(0,0, T)^\T$ and $(0,0,-T)^\T$. 

  Obviously, the optimal line has somehow to go through each cluster.
One possible candidate for the approximation line is
    simply the axis $\{(0,0,t) : t \in \BR\}$, whose distance to the
    given points is by construction $6$ -- for each cluster $3$.  Now,
    assume that the line goes through one given point, say
    $(1,0,T)^\T$.  If $T$ is very large, then we can neglect the slope
    of the line.  Only considering the distances within the cluster
    around $(0,0,T)^\T$, we notice that the distance increases from
    $3$ to approximately $2 \sqrt 3$.  Although the axis is maybe not
    the optimal line, the distance to the given points is smaller than
    for a line going through a data point.  Therefore, we can conclude
    that the optimal line has not to contain any given point.
  
  The same construction can be done for arbitrary subspaces
  of dimension
  $d  < n-1$.  For example, consider just
  the points
  \begin{equation*}
    (\cos(\phi),
    \sin(\phi) |
    \pm T \Vek e_k^\T)^\T
    \qquad\text{with}\qquad
    \phi \in \{0, \nicefrac{2\uppi}{3}, \nicefrac{4\uppi}{3}\},
    k = 1, \dots, d,
  \end{equation*}
  where $\Vek e_k$ is the $k$th unit vector.
  Using the same argumentation as above, the distance
    to the subspace $\mathrm{span} \{ \pm \Vek e_k : k=1,\dots, d\}$
    is smaller than to any $d$-dimensional subspace containing at
    least one data point.\qed
\end{Counterexample}

\section{Numerical examples}\label{sec:numerics}

In this section, we demonstrate the performance of r\pers{reaper} by numerical examples
implemented in MATLAB.

\subsection{(2,1)-norm versus Frobenius norm}\label{sec:numerical-examples}
\begin{figure}\centering
    \includegraphics
    {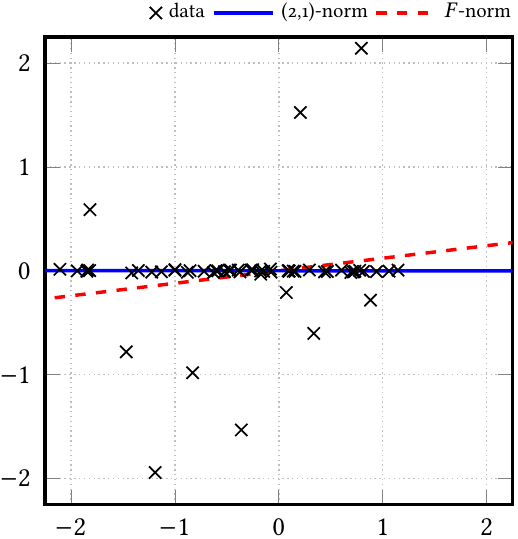}
  \caption{Performance of r\pers{reaper} with (2,1)-norm and Frobenius norm
   in the data fidelity term, respectively. The first one appears to be more robust against outliers.}
  \label{fig:syn1:dataset}
\end{figure}

This example with simple synthetic data will show that the (2,1)-norm in the data term is
more robust against outliers than the Frobenius norm
\begin{equation*}
  \biggl( \sum_{k=1}^N \pNormn{(\Mat P - \Mat I_n) \, \Vek x_k}_2^2
  \biggr)^\frac12
   = \pNormn{\Mat P \Mat X - \Mat X}_F.
\end{equation*}
For the Frobenius norm here abbreviated as $F$-norm, we have only to replace the projection 
onto $\mathfrak B_{2, \infty}$
with the projection to the Frobenius norm ball 
\begin{equation*}
  \proj_{\Set B_{F}}(\Mat Y)
  = \frac{\Mat Y}{[\pNormn{\Mat Y} - 1]_+
    + 1}.
\end{equation*}
We want to recover a line in the plane.  
Since this recovery problem is invariant under rotations,
we restrict ourselves to $\mathrm{span} \{ (1,0)^\T\}$.  The data are
generated randomly and consist of 50 points near the considered axis --
we added a small amount of noise in the second coordinate -- and of 10 outliers
located somewhere in the plane, see Figure~\ref{fig:syn1:dataset}.
Besides the data points, the recovered lines using r\pers{reaper} with
the (2,1)-norm (solid line) as data fidelity and the Frobenius norm
(dashed line) with parameters $d=1$ and $\alpha = 5$ are presented.
In this toy example, r\pers{reaper} yield nearly a perfect result
regardless of the outliers, and is in particular more robust than the
same model with the Frobenius norm.

\subsection{Nuclear norm and truncated hypercube constraints}
\label{sec:synthetic-subspace}

\begin{figure}\centering
  \subfloat[Representation of the dataset.]{
    \includegraphics
    {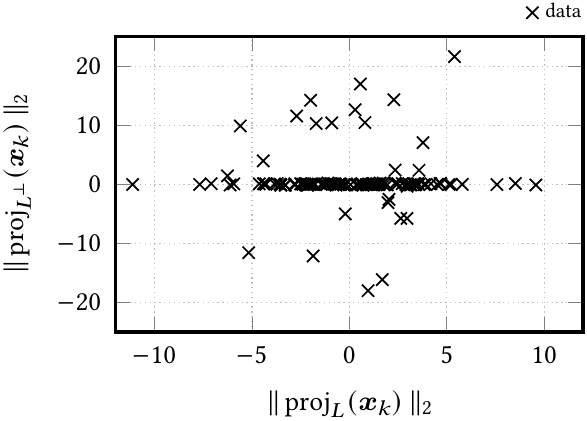}}
  \quad
  \subfloat[Evolution of the rank.]{
    \includegraphics
    {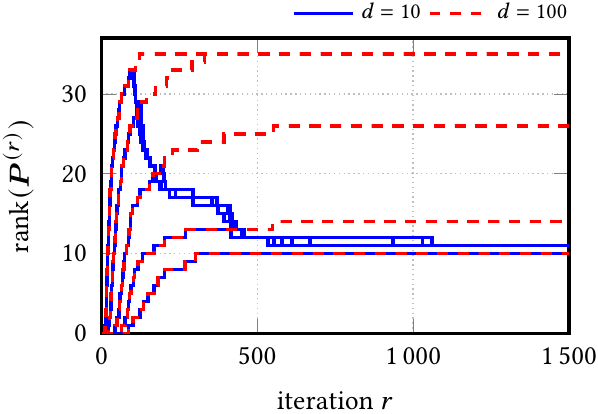}}
  \caption{Performance  of r\pers{reaper} for different upper dimension estimators
        $d=10,100$ and
        regularization parameters
        $\alpha = 2.5, 5, 10, 15, 20$ (top down for each $d$).}
  \label{fig:syn2}
\end{figure}

In this example we are interested how the
rank reduction is influenced by the nuclear norm and the projection to the
truncated hypercube.  In this synthetic  experiment, we approximate the given data
$\Vek x_k \in \BR^{100}$ by a 10-dimensional subspace.  The data is
again generated randomly, where 100 points lie near the subspace $L$
spanned by the first ten unit vectors and additional 25 outliers.  In
Figure~\ref{fig:syn2}.a, the dataset is represented by the distance to
the subspace $L$ and to the orthogonal complement $L^\perp$.

We apply r\pers{reaper} in
\thref{alg:matrix-free-reaper} with different parameter combination.
For the target dimension, we choose in our first experiment $d=10$, which is the wanted
dimension, and second one $d=100$, which does not truncate the unit hypercube
at all.  The influence of the regularization parameter $\alpha$ on the
rank of $\Mat P^{(r)}$ is shown in Figure~\ref{fig:syn2}.b, where the
lines from top to down correspond to the
regularization parameters $\alpha = 2.5, 5, 10 , 15, 20$.  
Since we start the
iteration with the zero-rank matrix $\Mat P^{(0)} \coloneqq \Mat 0$,
the first iterations for $d=10$ and $d=100$ coincides up to the point,
where the trace of $\Mat P^{(r)}$ exceeds the value 10.

Considering only the results for $d=100$ (solid lines), we see that
the nuclear norm reduces the rank of the iteration variable
$\Mat P^{(r)}$ with an increasing regularization parameter.
Further, the rank during the primal-dual algorithm is very sensitive
to the regularization parameter.  For $d=10$ (dashed lines), the
situation changes dramatically.  After the initial stages, the rank of
$\Mat P^{(r)}$ decreases nearly to the target dimension.  Since the
matrices $\Mat P^{(r)}$ are no orthogonal projections, rank and trace
do not conincide.  Due to this fact, the rank is not strictly bounded
by the maximal trace of the truncated hypercube.  Nevertheless, the
projection to the truncated hypercube significantly reduces the rank.

For an optimal rank evolution during the matrix-free primal-dual
method, the projection to the truncated hypercube by
\thref{alg:proj-trunc-cube} appers to be important.  Moreover, the
projection makes the rank evolution less sensitive to the
regularization parameter $\alpha$ so that a wider range of
regularization parameters can be applied without lossing the
computational benefits of the low rank.  Thus, the truncated hypercube
projection is an elementary key component of the algorithm.

\subsection{Face approximation} \label{sec:face-approx}

The idea to use the principle components of face images -- the
so-called eigenfaces -- for recognition, classification, and
reconstruction was considered in various paper and goes probably back to \cite{TP91}.
In this experiment, we adopt this idea to show that our matrix-free
\pers{reaper} can handle high-dimensional data.  Since the
computation of an optimal offset is non-trivial as discussed in Section
\ref{sec:incorporating-offset}, we choose just the geometric median.

\begin{figure}\centering
  \subfloat[Given full dataset of the experiment.]{
    \includegraphics[width=0.45\linewidth]
    {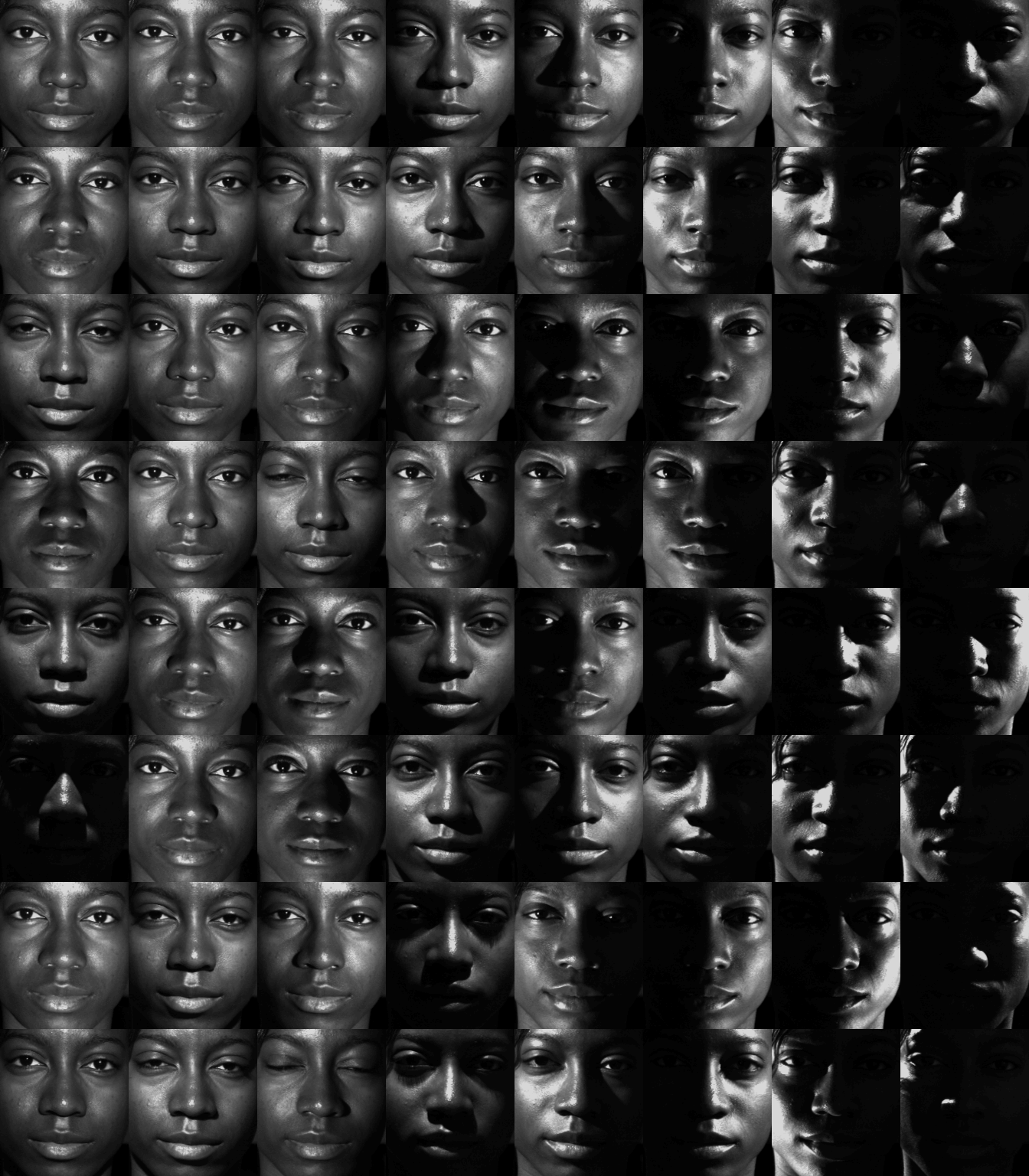}}
  \quad
  \subfloat[Projections to the determined subspace.]{
    \includegraphics[width=0.45\linewidth]
    {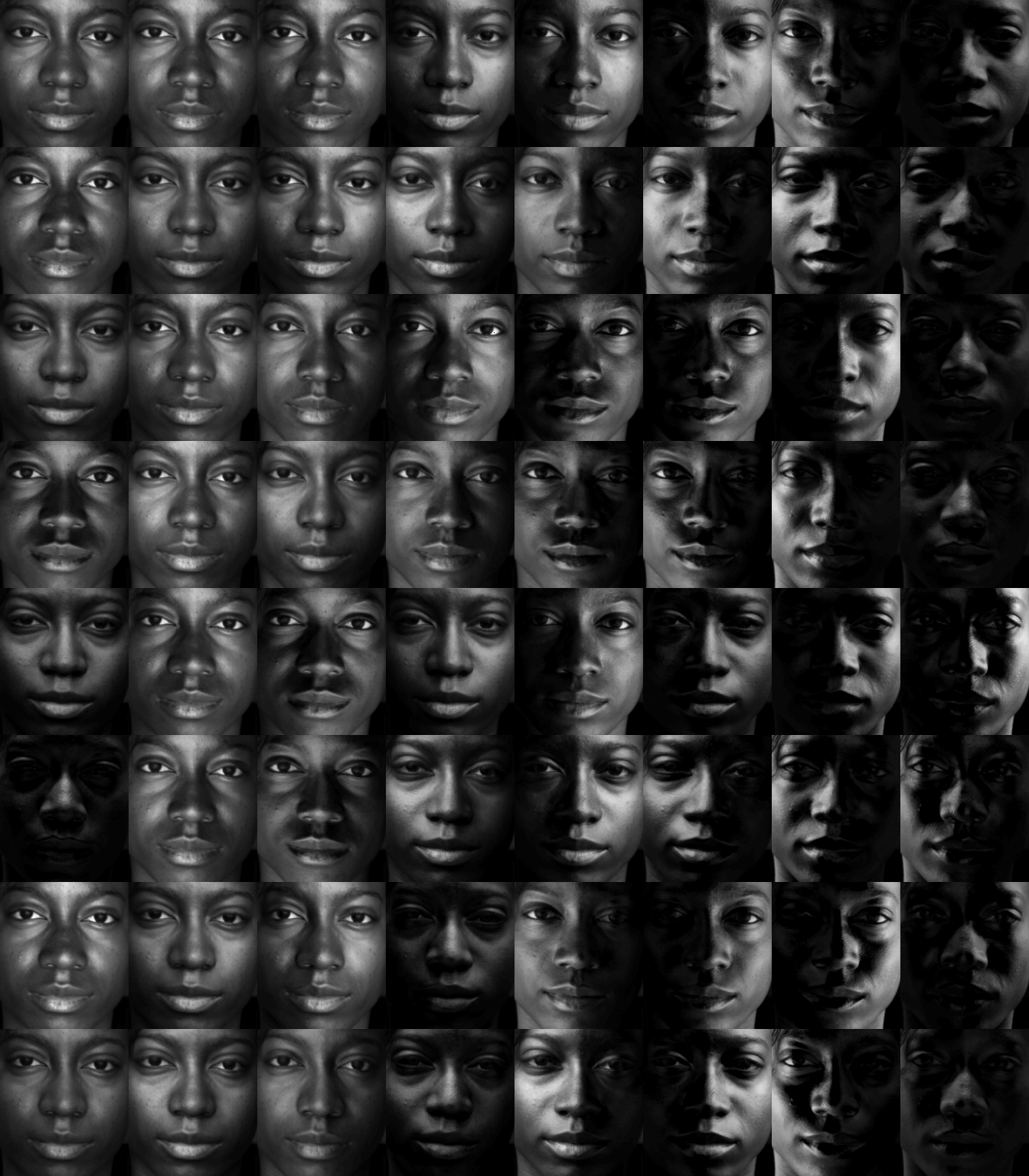}}
  \caption{The used images of the Extended Yale Face Dataset B and their projections
        onto the subspace determined by our matrix-free r\pers{reaper}.}
  \label{fig:ex1:dataset}
\end{figure}

\begin{figure}\centering
  \subfloat[Projection of a given data point.]{
    \includegraphics[width=0.2\linewidth]
    {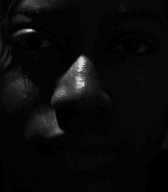}
    \quad
    \includegraphics[width=0.2\linewidth]
    {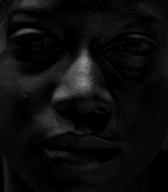}}
  \qquad
  \subfloat[Projection of an additional image.]{
    \includegraphics[width=0.2\linewidth]
    {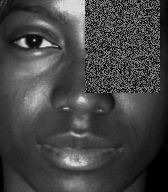}
    \quad
    \includegraphics[width=0.2\linewidth]
    {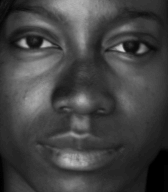}}
  \caption{Example of projections onto the recovered subspace.}
  \label{fig:ex1:detail}
\end{figure}

For the \textbf{first experiment}, we use the cropped \bq{Extended
  Yale Face Dataset B} \cite{GBK01, LHK05}.  The considered part of
the dataset consists in 64 images with 168$\times$192 pixels with
integer values between 0 and 255 of one face under different lighting
conditions, but with the same facial expression, see
Figure~\ref{fig:ex1:dataset}.a.
It is well-known that such images can be well approximated by a
subspace covering around five directions \cite{EHY95}.  In our
simulation, we set the maximal dimension to
$d=10$.  For the chosen regularization parameter
$\alpha = 2 \cdot 10^{4}$, our matrix-free r\pers{reaper} finds a
seven-dimensional subspace.  The projection of the original data to
this subspace is shown in Figure~\ref{fig:ex1:dataset}.b -- a higher
resolved example in Figure~\ref{fig:ex1:detail}.a.
An typical effect of the projection to the low-dimensional subspace is
that dark regions are lightened, shadows are removed, and reflections
at skin and eyes are cleared away.  The recovered subspace learned
from uncorrupted face images can be used to remove corrupted parts in
additional images as shown in Figure~\ref{fig:ex1:detail}.b.

\begin{figure}\centering
  \subfloat[Given full dataset of the experiment.]{
    \includegraphics[width=0.9\linewidth]
    {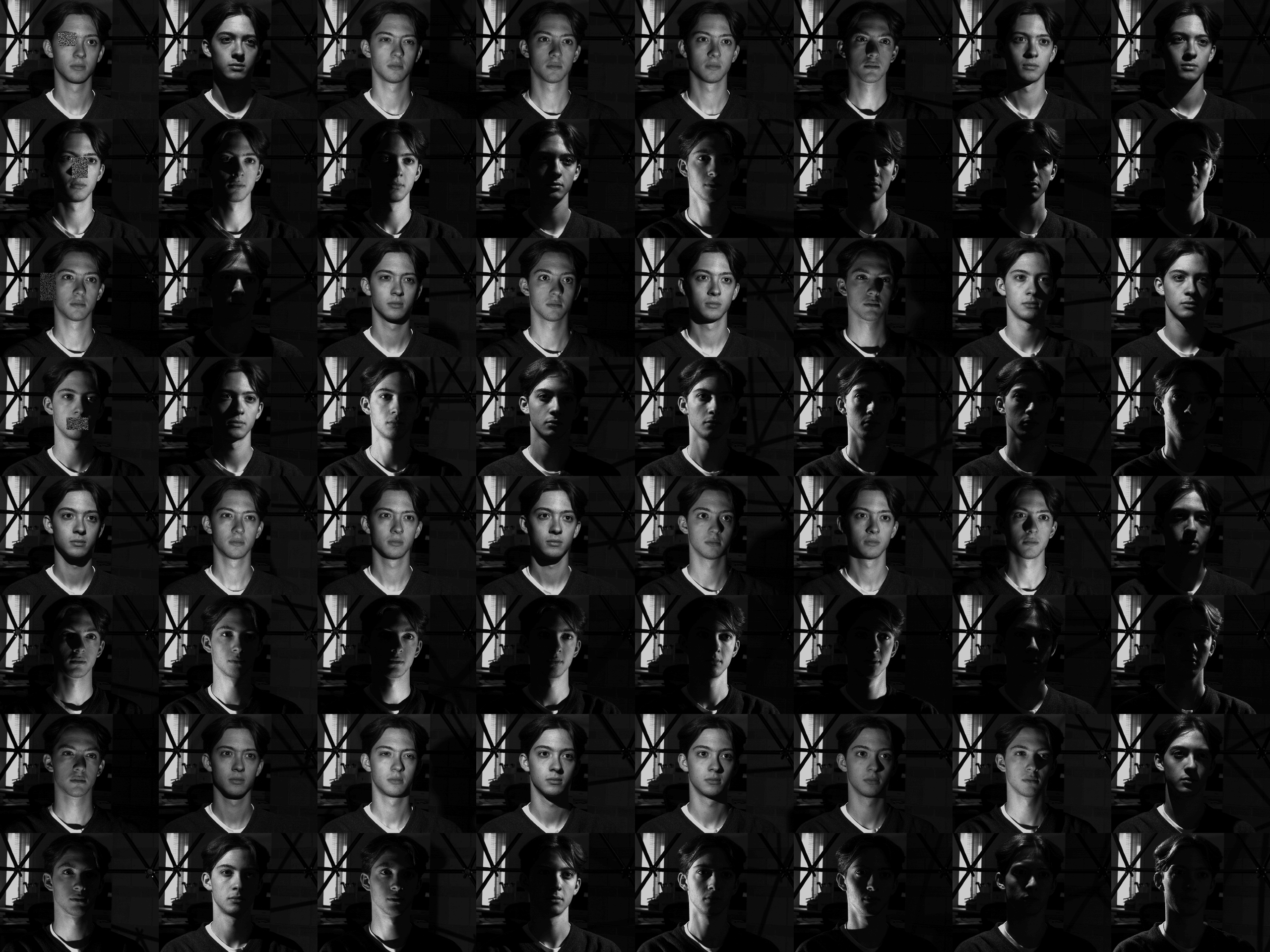}}
  \\[10pt]
  \subfloat[Projections to the determined subspace.]{
    \includegraphics[width=0.9\linewidth]
    {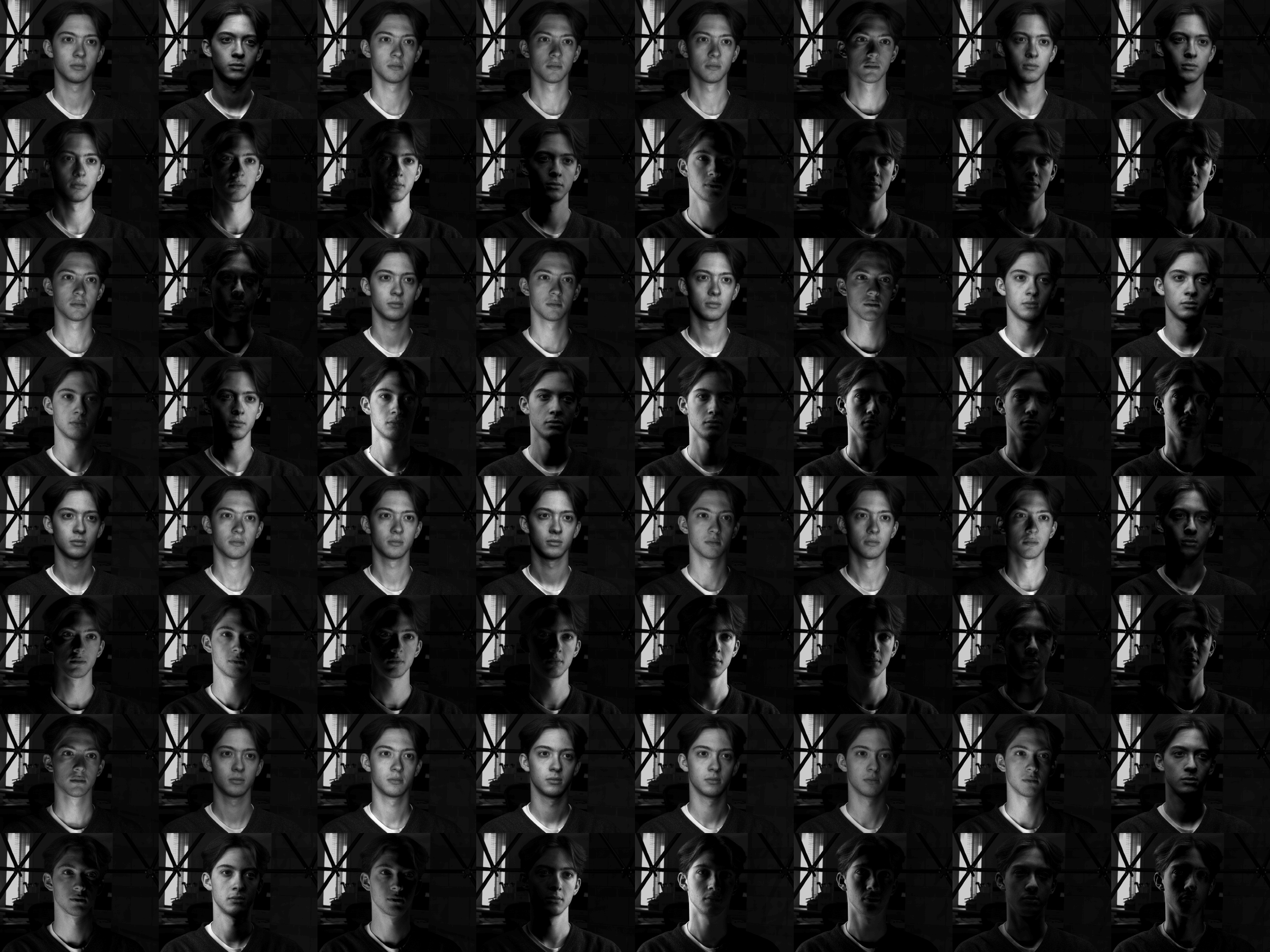}}
  \caption{The used images of the Extended Yale Face Dataset B and their projections
        onto the subspace determined by our matrix-free r\pers{reaper}.}
  \label{fig:ex3:dataset}
\end{figure}

In our \textbf{second experiment}, we consider images with a higher resolution.
The main motivation to develop a matrix-free algorithm have been to
handle such data.  We 
apply r\pers{reaper} to determine a five-dimensional subspace form the
full Extended Yale Face Dataset B.  
The used dataset is shown in
Figure~\ref{fig:ex3:dataset}, where each image has 640$\times$480
pixels.  Notice that an artefact has been placed in the first four
images covering the right eye, the nose, the right ear, and the mouth
respectively.

\begin{figure}\centering
    \includegraphics[width=0.4\linewidth]
    {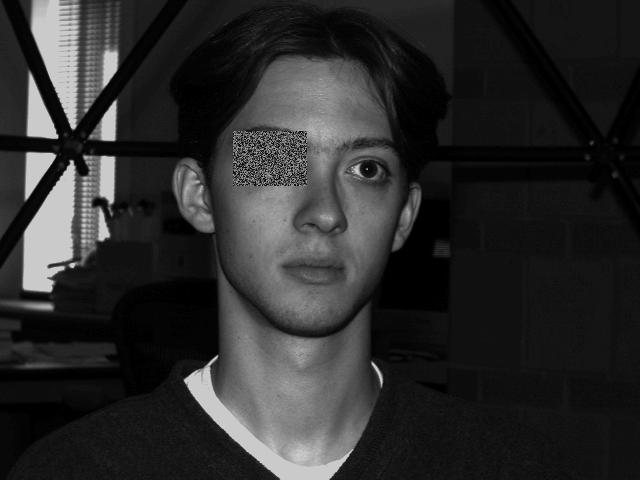}
    \quad
    \includegraphics[width=0.4\linewidth]
    {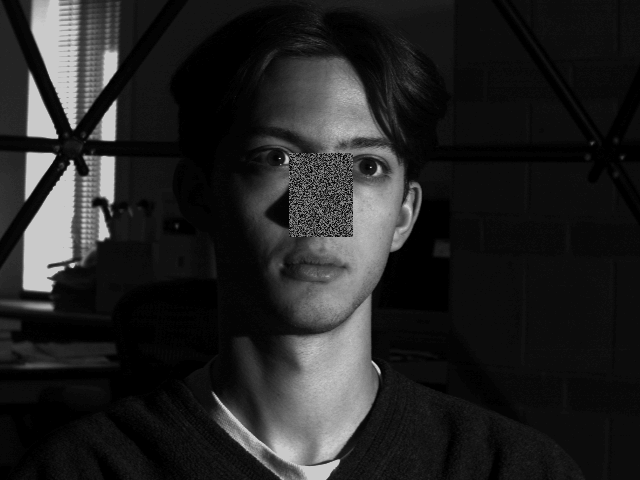}
    \\[20pt]
    \includegraphics[width=0.4\linewidth]
    {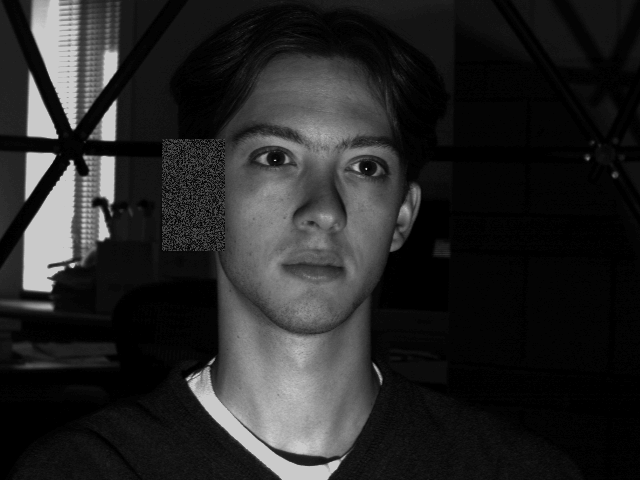}
    \quad
    \includegraphics[width=0.4\linewidth]
    {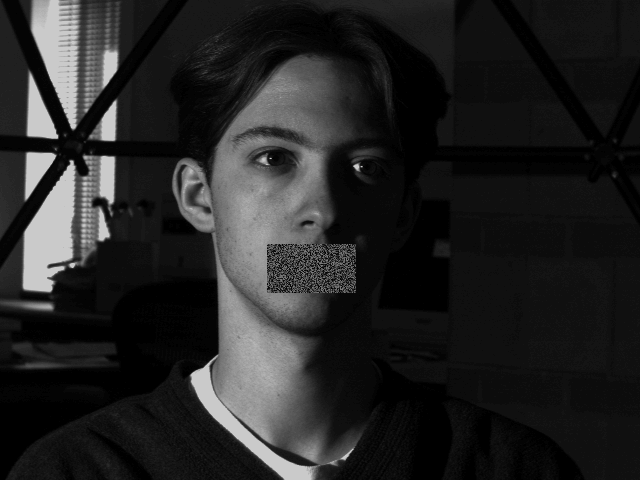}
  \caption{Corrupted images within the dataset.}
  \label{fig:ex3:arte}
\end{figure}

\begin{figure}\centering
    \includegraphics[width=0.4\linewidth]
    {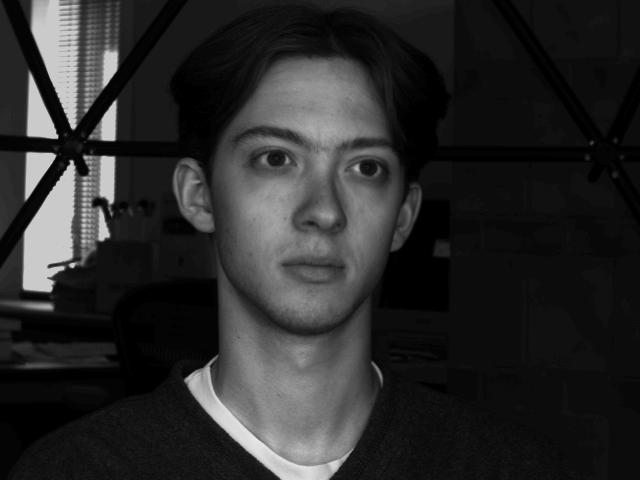}
    \quad
    \includegraphics[width=0.4\linewidth]
    {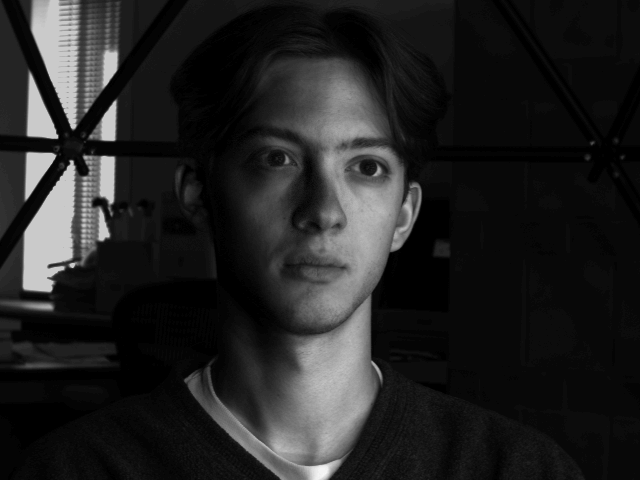}
    \\[20pt]
    \includegraphics[width=0.4\linewidth]
    {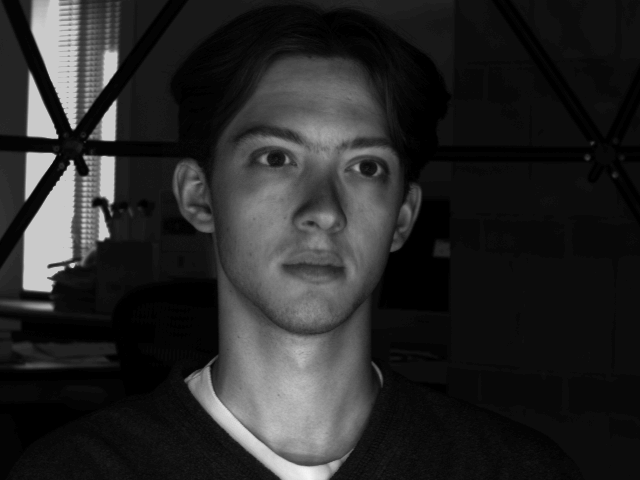}
    \quad
    \includegraphics[width=0.4\linewidth]
    {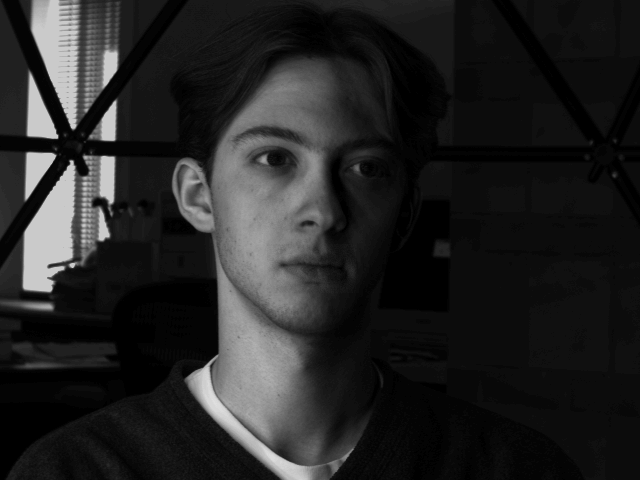}
  \caption{Restoration by projecting to the principle components.}
  \label{fig:ex3:detail}
\end{figure}

In order to remove the artifacts by unsupervised learning, we
approximate the full dataset including the artificial face images by a
five-dimensional subspace ($d=5$) using r\pers{reaper}, which should
be robust against the four outliers.  Projecting the first four images
to the recovered subspace, we removed the unwanted
artifacts, see Figure~\ref{fig:ex3:dataset}.b and
\ref{fig:ex3:arte}--\ref{fig:ex3:detail}.  

\begin{figure}\centering
  \subfloat[Rank.]{
    \includegraphics
    {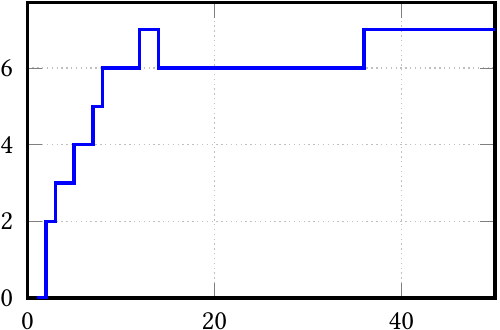}}
  \quad
  \subfloat[Objective.]{
    \includegraphics
    {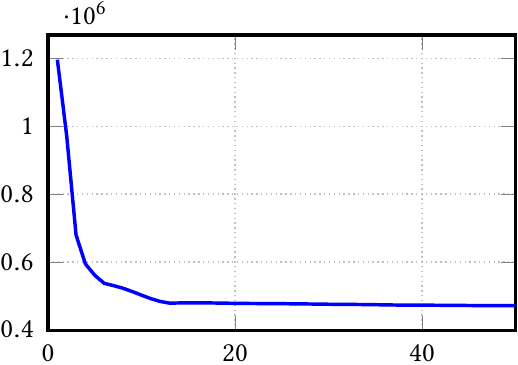}}
  \caption{Evolution of the rank and the objective during matrix-free
    r\pers{reaper}.}
  \label{fig:ex3:iter}
\end{figure}

Note that in this example the projection $\hat{\Mat \Pi}$ corresponds
to a 307\,200$\times$307\,200 matrix, which would require
703{.}125~GiB for double precision whereas the matrix-free
representation only requires around 16{.}407~MiB since the rank of the
primal variable is here bounded by seven, see
Figure~\ref{fig:ex3:iter}.  Further, we want to mention that the primal-dual minimization algorithm 
for r\pers{reaper} converges already after few iterations.

\section{Conclusion} \label{sec:conclusion}
Convex models are usually preferable over non-convex ones due to their
unique local minimum.  While robust PCA models that can handle high
dimensional data are usually nonconvex, a convex relaxation was
proposed by the \pers{reaper} model.  Relying on the projector
approach it is however not applicable for high dimensional data in its
original form.  To manage such data, we have combined primal-dual
optimization techniques from convex analysis with sparse factorization
techniques from the Lanczos algorithm.  Moreover, we extended the
model by penalizing the nuclear norm of the operator which has the
advantage that the dimension of the low dimensional subspace must not
be known in advance.  We addressed the problem of the bias in robust
PCA, but more research in this directions appears to be
necessary. Further other sparsity promoting norms then the nuclear norm
could be involved. Our method can be enlarged to 3D images as videos,
3D stacks of medical or material images, where tensor-free methods
will come into the play.  Finally, it may be interesting to couple PCA
ideas with approaches from deep learning to better understand the
structure of both.

\section*{Acknowledgment}
The authors want to thank G. Schneck for providing the idea of the proof for Theorem~\ref{thm:offset}.
 Funding by the German Research Foundation (DFG) with\-in the project
 STE 571/16-1 and by the Austrian Science Fund (FWF) within the
 project P28858 is gratefully acknowledged.

\end{document}